\documentclass[12pt]{article}
\RequirePackage{fix-cm}
\usepackage[T1]{fontenc}
\usepackage[utf8]{inputenc}
\usepackage{lmodern}
\usepackage{textcomp}
\usepackage{microtype}
\usepackage{authblk}

\usepackage[english]{babel}

\usepackage{amsmath}
\usepackage{amsfonts}
\usepackage{amssymb}
\usepackage{amsthm}
\usepackage{latexsym}
\usepackage{mathtools}
\usepackage{multirow}
\usepackage{array}
\usepackage{tabulary}

\usepackage{geometry} 
\usepackage[toc,page]{appendix}
\usepackage{lineno} 

\usepackage{xcolor}
\usepackage{float} 
\usepackage{placeins} 

\usepackage{csquotes}
\usepackage{url}
\usepackage{natbib}
\setcitestyle{authoryear,open={(},close={)}}
\newcommand{\textcite}{\cite}

\usepackage{graphicx}
\usepackage{setspace}

\usepackage{thmtools}

\declaretheorem[
	name=Theorem,
	numberwithin=section
	]{theorem}
\declaretheorem[
	name=Lemma,
	sibling=theorem,
	]{lemma}
\declaretheorem[
	name=Proposition,
	sibling=theorem,
	]{proposition}
\declaretheorem[
	name=Corollary,
	sibling=theorem,
	]{corollary}
\declaretheorem[
	name=Remark,
	sibling=theorem,
	style=remark
	]{remark}

\def\a{\alpha}
\newcommand{\be}{\beta}
\newcommand{\ep}{\epsilon}
\newcommand{\de}{\delta}
\newcommand{\ga}{\gamma}

\newcommand{\la}{\lambda}
\newcommand{\rh}{\rho}
\newcommand{\si}{\sigma}
\newcommand{\ze}{\zeta}

\newcommand{\De}{\Delta}
\newcommand{\Th}{\Theta}

\newcommand\Pinf{P^\infty}
\newcommand\Pexn{P^{(n)}}
\newcommand\Pn{P^{(n)}}
\newcommand\PnFth{P_{\rm{F3}}^{(n)}}
\newcommand\PinfFth{P_{\rm{F3}}^\infty}
\newcommand\PinfFL{P_{\rm{FL}}^\infty}
\newcommand\PnFL{P_{\rm{FL}}^{(n)}}
\newcommand\PinfPoi{P_{\rm{Poi}}^\infty}
\newcommand\PnPoi{P_{\rm{Poi}}^{(n)}}
\newcommand\Pinfphi{P_{\varphi}^\infty}
\newcommand\Pnphi{P_{\varphi}^{(n)}}
\newcommand\PinfBin{P_{\rm{Bin}}^\infty}

\newcommand\PnNB{P_{\rm{NB}}^{(n)}}
\newcommand\PinfNB{P_{\rm{NB}}^\infty}
\newcommand\PinfGP{P_{\rm{GP}}^\infty}
\newcommand\PnGP{P_{\rm{GP}}^{(n)}}
\newcommand\Sinf{S^\infty}
\newcommand\Sinfphi{S_{\varphi}^\infty}
\newcommand\Sn{S^{(n)}}
\newcommand\Snphi{S_{\varphi}^{(n)}}
\newcommand\SinfFL{S_{\rm{FL}}^\infty}
\newcommand\SnFL{S_{\rm{FL}}^{(n)}}
\newcommand\SnPoi{S_{\rm{Poi}}^{(n)}}
\newcommand\SinfPoi{S_{\rm{Poi}}^\infty}

\newcommand\SnNB{S_{\rm{NB}}^{(n)}}
\newcommand\SinfBin{S_{\rm{Bin}}^\infty}
\newcommand\SinfNB{S_{\rm{NB}}^\infty}
\newcommand\SnFth{S_{\rm{F3}}^{(n)}}

\newcommand\SinfGP{S_{\rm{GP}}^\infty}
\newcommand\SnGP{S_{\rm{GP}}^{(n)}}
\newcommand\UQ{U^{\rm{Q}}}
\newcommand\LQ{L^{\rm{Q}}}
\newcommand\UDN{U^{\rm{DN}}}
\newcommand\gaFth{\gamma_{\rm{F3}}}
\newcommand\gaFL{\gamma_{\rm{FL}}}
\newcommand\gaPoi{\gamma_{\rm{Poi}}}
\newcommand\gaBin{\gamma_{\rm{Bin}}}
\newcommand\gaphi{\gamma_{\varphi}}
\newcommand\gaNB{\gamma_{\rm{NB}}}
\newcommand\gaGP{\gamma_{\rm{GP}}}
\newcommand\varphiFL{\varphi_{\rm{FL}}}
\newcommand\varphiPoi{\varphi_{\rm{Poi}}}
\newcommand\varphiFth{\varphi_{\rm{F3}}}
\newcommand\varphiBin{\varphi_{\rm{Bin}}}
\newcommand\varphiNB{\varphi_{\rm{NB}}}
\newcommand\varphiGP{\varphi_{\rm{GP}}}
\newcommand\fP{f_{\rm{Poi}}}

\newcommand\fNB{f_{\rm{NB}}}
\newcommand\gNB{g_{\rm{NB}}}
\newcommand\fFth{f_{\rm{F3}}}
\newcommand\fGP{f_{\rm{GP}}}
\newcommand\fBinFL{f_{\rm{Bin}}}
\newcommand\ffBinFL{\hat f_{\rm{Bin}}}
\newcommand\rFth{\rh_{\rm{F3}}}
\newcommand\pFth{\pi_{\rm{F3}}}
\newcommand\rBin{\rh_{\rm{Bin}}}
\newcommand\pBin{\pi_{\rm{Bin}}}
\newcommand\rNB{\rh_{\rm{NB}}}
\newcommand\pNB{\pi_{\rm{NB}}}

\newcommand\mFL{m_{\rm{FL}}}
\newcommand\mBin{m_{\rm{Bin}}}
\newcommand\mNB{m_{\rm{NB}}}
\newcommand\mFth{m_{\rm{F3}}}
\newcommand\mphi{m_{\varphi}}

\newcommand\pnr{p_0^{(r)}}
\newcommand\pnga{p_0^{(\ga)}}
\newcommand\pnplus{p_0^{(+)}}

\newcommand\xit{{\tilde\xi}}

\newcommand\zet{{\tilde\ze}}
\newcommand\yx{y(x)}
\newcommand{\EV}[1]{\operatorname{E}#1}
\newcommand{\Var}[1]{\operatorname{Var}#1}

\newcommand{\Prob}[1]{\operatorname{Prob}#1}

\newcommand{\ProdLog}[1]{\operatorname{ProductLog}#1}
\defcitealias{Agresti1974}{Agresti's (1974)}
\defcitealias{Athreya1992}{Athreya's (1992)}
\defcitealias{Haldane1927}{Haldane's (1927)}
\defcitealias{Hill1982a}{Hill's (1982)}
\defcitealias{Hoellinger2019}{H\"ollinger's (2019)}
\defcitealias{Lande1976}{Lande's (1976) equation}
\defcitealias{Mathematica}{{\it Mathematica}}
\defcitealias{Quine1976}{Quine's (1976)}
\defcitealias{Pollak1971}{Pollak's (1971)}
\numberwithin{equation}{section}
\numberwithin{figure}{section}
\numberwithin{table}{section}
\begin{document}
\normalbaselineskip18pt
\baselineskip18pt
\global\hoffset=-10truemm
\global\voffset=-10truemm
\allowdisplaybreaks[2]
\setcounter{page}{1}
%
%
\title{\bf Bounds for survival probabilities in \\ supercritical Galton-Watson processes and \\ applications to population genetics \\ 
\vphantom{T}
}
\author{\bf Reinhard B\"urger}

\affil{Faculty of Mathematics \\ University of Vienna \\ Oskar-Morgenstern-Platz 1 \\ 1090 Vienna, Austria  \\
\vspace{0.5cm} 
Email: reinhard.buerger@univie.ac.at \\
ORCID: 0000-0003-2174-4968
}
\date{}

\maketitle

\vspace{4cm}
\thispagestyle{empty}


\vspace{0.5cm}

\vspace{0.5cm}
\noindent{\bf Additional materials} A \emph{Mathematica} notebook containing the elements to check the algebraically demanding parts of the proofs and to produce figures and tables is provided as Supplementary Material. It can also be obtained directly from the author.

\newpage
\begin{abstract}
Population genetic processes, such as the adaptation of a quantitative trait to directional selection, may occur on longer time scales than the sweep of a single advantageous mutation. To study such processes in finite populations, approximations for the time course of the distribution of a beneficial mutation were derived previously by branching process methods. The application to the evolution of a quantitative trait requires bounds for the probability of survival $S^{(n)}$ up to generation $n$ of a single beneficial mutation. Here, we present a method to obtain a simple, analytically explicit, either upper or lower, bound for $S^{(n)}$ in a supercritical Galton-Watson process. We prove the existence of an upper bound for offspring distributions including Poisson, binomial, and negative binomial. They are constructed by bounding the given generating function, $\varphi$, by a fractional linear one that has the same survival probability $S^\infty$ and yields the same rate of convergence of $S^{(n)}$ to $S^\infty$ as $\varphi$. For distributions with at most three offspring, we characterize when this method yields an upper bound, a lower bound, or only an approximation. 
Because for many distributions it is difficult to get a handle on $S^\infty$, we derive an approximation by series expansion in $s$, where $s$ is the selective advantage of the mutant. We briefly review well-known asymptotic results that generalize Haldane's approximation $2s$ for $S^\infty$, as well as less well-known results on sharp bounds for $S^\infty$.
We apply them to explore when bounds for $S^{(n)}$ exist for a family of generalized Poisson distributions. Numerical results demonstrate the accuracy of our and of previously derived bounds for $S^\infty$ and $S^{(n)}$. Finally, we treat an application of these results to determine the response of a quantitative trait to prolonged directional selection.
\end{abstract}
\section{Introduction}
Galton-Watson branching processes were used early in the history of population genetics to approximate the fixation probability of a single advantageous mutation in a finite population \citep{Fisher1922,Haldane1927}. In particular, Haldane showed that if the offspring distribution of the mutant is approximately Poisson with mean $m=1+s$, the fixation probability can be approximated by $2s$ provided the selective advantage $s$ is sufficiently small. More recent work led to considerable generalizations of Haldane's approximation and is discussed in Sect.~\ref{subsec:bounds_series}.

Essentially parallel to the first applications of Galton-Watson processes the method of diffusion approximation was introduced and developed \citep{Fisher1922,Wright1931,Kolmogorov1931}.  Whereas this method is a powerful tool to quantify fixation probabilities, stationary distributions, or the distribution of the time to fixation of an allele \citep{Kimura1964,Ewens2004}, it is less well suited to obtain analytically explicit expressions for the time-dependence of allele frequencies under selection. Approximations have been derived, essentially for statistical purposes, but they are semi-explicit and complex \citep[e.g.][]{Steinrucken_etal2013}. 

The time course of the frequency distribution of a new favorable mutant has been approximated by explicit formulas derived with the help of branching-process theory \citep{DesaiFisher2007, UeckerHermisson2011, MartinLambert2015, GB2024}. Application of these results to population genetic processes that occur on longer time scales than the sweep of a single mutation, such as the evolutionary response of a quantitative trait to directional selection, require bounds for the probability of survival $\Sn$  of a single new mutant up to generation $n$ \citep{GB2024}. 
The response of the mean of a trait is determined by the variance contributed by every mutation that is favored by selection and spreads. The total variance contributed by a single mutation while present in the population is given by an integral, whose integrand depends, among others, on $\Sn$. With accurate and analytically simple bounds on $\Sn$, this integral can be approximated and the error estimated. This procedures also requires estimates on the time $T(\ep)$ needed for $\Sn$ to fall below $(1+\ep)\Sinf$, $\Sinf$ the (ultimate) survival probability. Under appropriate scaling assumptions on the strength $s$ of selection and the population size $N$, the time $T(\ep)$ is short compared with the time while the mutant is sweeping to fixation. Essentially, it can be shown that the variance contributed by mutations that are lost is negligible. Above $T(\ep)$, $\Sn$ can be approximated by the constant term $\Sinf$ which greatly simplifies the integral (for a detailed description, see Sect.~\ref{subsec:spread_favorable_mutant}).

\citet{GB2024} imposed the assumption that the offspring distribution is such that $\Sn$ can be bounded above by the explicitly available expression for an appropriately chosen modified geometric distribution (which has a fractional linear generating function; see Sect.~\ref{sec:FL}). As explained below, here we present a general method to derive such bounds and prove its applicability for some well-known families of offspring distributions.

Motivated by these considerations, the main goal of this paper is the derivation of sharp, explicit, and analytically tractable upper bounds for the probability of survival $\Sn$ up to generation $n$ in supercritical Galton-Watson processes. We adopt the principal method pioneered by \citet{Seneta1967} (and attributed by him to P.A.P.~Moran) of using probability generating functions (pgfs) of fractional linear type to bound a given pgf $\varphi$. \citet{Seneta1967} and \citet{Agresti1974} used it to derive simple bounds for the extinction time distribution of subcritical or critical branching processes, originating from specific offspring distributions. For the Poisson distribution, Agresti derived best possible bounds of fractional linear type and indicated how to derive bounds for the supercritical case by exploiting a duality relation between subcritical and supercritical processes. For the supercritical case, these bounds are no longer pgfs (see Sect.~\ref{sec:alternatives}). A different method to obtain bounds for the extinction probabilities $\Pn=1-\Sn$ for a given pgf was developed by \citet{Pollak1971}. It is based on series expansion of the pgf and is applicable to sub- and supercritical processes (see Sect.~\ref{sec:alternatives}).

We use a direct method for the supercritical case that is based on proper generating functions and requires that the prospective bounding fractional linear pgf has the same extinction probability $\Pinfphi$ and the same slope $\gaphi=\varphi'(\Pinfphi)$ as the given pgf $\varphi$ (see Sect.~\ref{sec:basic_result}). This method can be applied to pgfs other than the Poisson distribution, even if no analytical expression for $\Pinfphi$ is available, as in the case of a binomial distribution when a re-parameterization in terms of $\PinfBin$ and the number of trials is possible. 

Using this method, we prove that simple, explicit upper bounds obtained from fractional linear distributions, denoted $\SnFL$, do exist for Poisson, binomial, and negative binomial distributions (Sects.~\ref{sec:Poisson} -- \ref{sec:NB}). For distributions with at most three offspring, $\SnFL$ can yield an upper bound (in most of the parameter space), a lower bound, or $\SnFL$ may switch from a lower to an upper bound at some generation $n$; a full characterization is obtained in Sect.~\ref{sec:F3}. Except for the Poisson distribution, where the proof is simple enough to provide insight, the proofs are relegated to the Appendix.

For most distributions our method is difficult to apply because their pgfs are too complicated to be handled analytically and already $\Pinf$ is difficult to access. Interestingly, there exists a branch of research that seems to be completely disconnected from the literature related to Haldane's approximation. 
\citet{Quine1976}, \citet{DaleyNarayan1980}, \citet{From2007}, and others derived explicit and very accurate upper and lower bounds for the eventual extinction probability in Galton-Watson processes with offspring distributions having finite second or third moment. They are outlined in Sects.~\ref{sec:Quine_bound} and \ref{sec:DN_bound};  for an extensive review consult From's paper. In Sect.~\ref{subsec:bounds_series}, we devise a method for general pgfs $\varphi$ to deduce series expansions of $\Pinfphi$ and of $\gaphi$ in terms of $s$, where $m=1+s>1$ and $s$ is small. The bounds of \citet{Quine1976} and \citet{DaleyNarayan1980} are shown to have an error of order $O(s^3)$. In this context we also briefly review recent, far reaching generalizations of Haldane's approximation for the fixation probability in finite populations. In Sect.~\ref{subsec:fix_prob_WF}, we review bounds and approximations obtained previously by diffusion-approximation methods for the Wright-Fisher model.

Among others, we apply our series expansions to derive analytically explicit, at least approximate, bounds for $\Snphi$ for a family of generalized Poisson distributions which otherwise is prohibitively difficult to tackle (Sect.~\ref{subsec:GP}). Even for cases that can be treated fully analytically, such as the Poisson or binomial distribution, these expansions yield valuable additional insights (Sects.~\ref{subsec:Poisson_bounds}, \ref{subsec:Bin_bounds}, \ref{subsec:NB_bounds}).  For the generalized Poisson distribution, $\SnFL$ may yield an upper bound for the true $\SnGP$ (if the variance is not much higher than the mean), a lower bound (if the variance is much higher than the mean), or switch from an upper to a lower bound at some $n$ (Sect.~\ref{subsec:GP}).
 
For every pgf $\varphi$ with finite variance, the sequence $\SnFL$ that we construct, whether it is an exact bound or an approximation based on series expansion, has the property that it converges to the given survival probability $\Sinfphi$ at the correct asymptotic rate $\gaphi$. The accuracy of the resulting convergence times $T_\varphi(\ep)$ and relative errors of the bounds for $\Snphi$ are explored in Sects.~\ref{subsect:Tep} and \ref{sec:appl_ext_times}, respectively.
Our central population genetics application is treated in Sect.~\ref{subsec:spread_favorable_mutant}.

\section{Definitions and preliminaries}\label{Sec:basics}
\subsection{Basic notation and assumptions}
We consider a Galton-Watson process $\{Z_n\}$, where $Z_n=\sum_{j=1}^{Z_{n-1}}\xi_j$, $Z_0=1$, and $\xi_j$ denotes the (random) number of offspring of individual $j$ in generation $n-1$. Thus, $Z_n$ counts the number of descendants of the mutant that emerged in generation 0. We assume that the $\xi_j$ are mutually independent, identically distributed random variables, independent of $n$, and have at least three finite moments, where 
\begin{linenomath}\begin{equation}\label{def_mv}
	\EV [\xi_j] = m>1 \;\text{ and }\; \Var[\xi_j] = \si^2 > 0\,.
\end{equation}\end{linenomath}
Therefore, the process $\{Z_n\}$ is supercritical. We denote the probability generating function (pgf) of $\xi_j$ (hence of $Z_1$) by $\varphi$, and the probability of having $k$ offspring by $p_k=P(\xi_j=k)$. Then $\varphi(x) = \sum_{k=0}^\infty p_k x^k$. To avoid trivialities, we always assume $p_0>0$ and $p_0+p_1<1$ (whence $\si^2>0$). 

As is well known \citep[e.g.][]{AthreyaNei1972}, the extinction probability $\Pinfphi$ of this process is the unique value $x\in(0,1)$ satisfying $\varphi(x)=x$. We denote by
\begin{linenomath}\begin{equation}\label{def_Pexn}
	\Pnphi=\Prob[Z_n=0|Z_0=1]
\end{equation}\end{linenomath}
the probability that the mutant is extinct by generation $n$. It satisfies $\Pnphi=\varphi^{(n)}(0)$, where $\varphi^{(n)}$ is the $n$th iterate of $\varphi$. Our assumptions imply $0<\Pnphi<\Pinfphi <1$ and $\lim_{n\to\infty} \Pnphi=\Pinfphi$. Often it will be convenient to formulate results in terms of  the corresponding survival probabilities:
\begin{linenomath}\begin{equation}\label{def_Sn}
	\Sinfphi = 1-\Pinfphi \; \text{ and }\; \Snphi=1-\Pnphi\,.
\end{equation}\end{linenomath}
We will use subscripts, such as $\varphiPoi$ and $\PinfPoi$, to refer to specific offspring distributions.

\subsection{Seneta's method of bounding the extinction probabilities $\Pnphi$}\label{Seneta's method}
\citet{Seneta1967} showed the following result (which does not require $m>1$):
Let $\varphi_L$, $\varphi_U$, and $\varphi$ be pgfs such that
\begin{linenomath}\begin{equation}
	\varphi_L(x) \le \varphi(x) \le \varphi_U(x) \;\text{for every }\; x\in[0,1] \,.
\end{equation}\end{linenomath}
Then 
\begin{linenomath}\begin{equation}
	\varphi_L^{(n)}(x) \le \varphi^{(n)}(x) \le \varphi_U^{(n)}(x) \;\text{ for every }\; x\in[0,1] \text{ and  every } n\ge1\,.
\end{equation}\end{linenomath}
In particular, the probability of extinction by generation $n$, $P_\varphi^{(n)}$, satisfies
\begin{linenomath}\begin{equation}\label{Senetas_inequ}
	\varphi_L^{(n)}(0) \le P_\varphi^{(n)} \le \varphi_U^{(n)}(0) \;\text{ for every }\; n\ge1\,.
\end{equation}\end{linenomath}

If the Galton-Watson process generated by $\varphi$ is supercritical, it is natural to bound it by supercritical processes. Because $\Pnphi$ is monotone increasing and converges to $\Pinfphi<1$, the following variant of this results is valid:
If
\begin{linenomath}\begin{equation}\label{Senetas_inequ_assump+}
	\varphi_L(x) \le \varphi(x) \le \varphi_U(x) \;\text{for every }\; x\in[0,\Pinfphi] \,,
\end{equation}\end{linenomath}
then \eqref{Senetas_inequ} holds. 
As in some previous work \citep[e.g.][]{Seneta1967,Agresti1974}, we will use fractional linear pgfs as bounds because they have the property that the $n$th iterates $\varphi^{(n)}$ can be calculated explicitly.

\subsection{Fractional linear generating functions}\label{sec:FL}
The modified geometric, or fractional linear, distribution is defined by
\begin{equation}\label{frac_lin}
	p^{\rm{(FL)}}_0 = \rh \;\text{ and }\; p^{\rm{(FL)}}_k = (1-\rh)(1-\pi)\pi^{k-1} \; \text{ if }\; k\ge1 \,,
\end{equation}
where $0<\rh<1$ and $0<\pi<1$ (e.g.~\citealt[pp.\ 6-7]{AthreyaNei1972}; \citealt[p.\ 16]{Haccou2005}). The name fractional linear derives from the fact that its pgf is
\begin{equation}\label{frac_lin_gen}
	\varphiFL(x;\pi,\rh) = \frac{\rh + x(1-\pi-\rh)}{1 - x \pi}\,,
\end{equation}
hence fractional linear. With $\rh=1-\pi$, the geometric distribution is recovered. It is straightforward to show that every fractional linear pgf generates a modified geometric distribution. We omit the dependence of $\varphiFL$ on $\pi$ and $\rh$ if no confusion can occur.

Mean and variance of $\{p^{\rm{(FL)}}_k\}$ are
\begin{equation}
	\mFL = \frac{1-\rh}{1-\pi} \text{ and } \si^2_{\rm{FL}} = \frac{(1-\rh)(\pi+\rh)}{(1-\pi)^2}\,.
\end{equation}
Therefore $\mFL>1$ if and only if $0<\rh<\pi<1$, and $\mFL>\si^2_{\rm{FL}}$ if and only if $2\pi+\rh < 1$, which implies $\pi<\tfrac12$.  

If $m\neq1$, and after rearrangement of the parameterization in \citet[p.~7]{AthreyaNei1972}, the $n$-times iterated pgf is again fractional linear and has parameters
\begin{equation}\label{pnrn}
	\pi_n = \frac{\pi(1-\mFL^{-n})}{\pi-\rh \mFL^{-n}} \text{ and } \rh_n = \frac{\rh(1-\mFL^{-n})}{\pi-\rh \mFL^{-n}}\,. 
\end{equation}

Now assume $\mFL>1$, i.e., $\rh<\pi$. Then the probability of extinction by generation $n$ is $\PnFL = \varphiFL^{(n)}(0) = \rh_n$ and the (ultimate) extinction probability is 
\begin{equation}\label{PinfFL}
	\PinfFL = \frac{\rh}{\pi}\,.
\end{equation}
By simple algebra we arrive at
\begin{linenomath}\begin{equation}\label{PnFL}
	\PnFL = \frac{\PinfFL(1-\mFL^{-n})}{1-\mFL^{-n}\PinfFL}
\end{equation}\end{linenomath}
and
\begin{linenomath}\begin{equation}\label{SnFL}
	\SnFL = 1-\PnFL = \frac{\SinfFL}{1- \mFL^{-n}(1-\SinfFL)} \,. 
\end{equation}\end{linenomath}
Because it will be important in subsequent sections, we note that
\begin{linenomath}\begin{equation}\label{gaFL=1/mFL}
	\gaFL := \varphiFL'(\PinfFL) = \mFL^{-1} \,. 
\end{equation}\end{linenomath}

Equation \eqref{SnFL} allows to compute the time needed for the probability of survival up to generation $n$, $\SnFL$, to decline to $(1+\ep)\SinfFL$.  For $\ep>0$ (not necessarily small) we define $T_{\rm FL}(\ep)$ as the (positive) solution $T$ of
\begin{linenomath}
\begin{equation}\label{Tep_definition_FL}
	S_{\rm FL}^{(T)} = (1+\ep)\SinfFL\,.
\end{equation}
\end{linenomath}
With the help of \eqref{SnFL}, this time is 
\begin{linenomath}\begin{equation}\label{Tep_FL}
	T_{\rm FL}(\ep) = \frac{\ln\Bigl(\bigl(1+\frac{1}{\ep} \bigr) \PinfFL \Bigr)}{\ln \mFL}  \,.
\end{equation}\end{linenomath}
Of course, the first generation in the associated GW-process that satisfies $S_{\rm FL}^{(T)} \le (1+\ep)\SinfFL$ is the least integer greater than or equal $T_{\rm FL}(\ep)$.

\section{The basic result and alternative methods for deriving bounds for $\Snphi$}\label{sec:General_approach}
First we derive our basic result and simple consequences. Then we discuss alternative approaches.

\subsection{Basic result}\label{sec:basic_result}
For a given pgf $\varphi$, we are primarily interested in lower bounds for $\Pnphi$, and upper bounds for $\Snphi$. It is well known that $\Pnphi$ converges to $\Pinfphi$ at the geometric rate 
\begin{linenomath}\begin{equation}\label{def_gamma}
	\gaphi := \varphi'(\Pinfphi)
\end{equation}\end{linenomath}
(Athreya and Nei 1972, Sect.~1.11). By Seneta's inequalities \eqref{Senetas_inequ}, we can obtain a lower bound for $\Pnphi$ that converges to $\Pinfphi$ at the correct rate $\gaphi$, if we can choose a fractional linear pgf $\varphiFL$ such that 
\begin{linenomath}\begin{equation}\label{cond_Pinf_gainf}
	\varphiFL(\Pinfphi)=\Pinfphi \;\text{ and }\; \gaFL=\varphiFL'(\Pinfphi)=\gaphi
\end{equation}\end{linenomath}
and $\varphiFL(x;\pi,\rh)\le \varphi(x)$ for every $x\in[0,\Pinfphi]$.  
Indeed, a straightforward calculation shows that for given $0<a_1<1$ and $0<a_2<1$, there is always a unique solution $(\pi,\rh)$ of the system
\begin{linenomath}\begin{equation}\label{pr_ab1}
	\varphiFL(a_1;\pi,\rh)=a_1 \text{ and } \varphiFL'(a_1;\pi,\rh)=a_2\,.
\end{equation}\end{linenomath}
 It is given by
\begin{linenomath}\begin{equation}\label{pr_solution_gam}
	\pi = \frac{1-a_2}{1- a_1a_2} \text{ and } \rh = a_1 \pi
\end{equation}\end{linenomath}
and satisfies $0<\rh<\pi<1$. 

With $a_1=\Pinfphi$, $a_2=\gaphi$, and the resulting values $\pi$ and $\rh$, eq.~\eqref{PnFL} informs us that for the resulting  fractional linear offspring distribution the probability of extinction by generation $n$ is
\begin{linenomath}\begin{equation}\label{PnFL_varphi}
	\PnFL =  \varphiFL^{(n)}(0) =  \frac{\Pinfphi(1-\gaphi^n)}{1-\gaphi^n\Pinfphi} \,,
\end{equation}\end{linenomath}
where we used \eqref{gaFL=1/mFL}.
Together with the left inequality in \eqref{Senetas_inequ}, these considerations yield the following basic result.

\begin{proposition}\label{prop:bound_Psurv_gamma}
Let $\varphi(x)$ be a pgf satisfying our general assumptions stated in Section \ref{Sec:basics}, so that $m>1$ and $0<\Pinfphi<1$. Let $\varphiFL(x;\pi_\varphi,\rh_\varphi)$ denote the uniquely determined fractional linear pgf that satisfies \eqref{cond_Pinf_gainf}. 
If 
\begin{linenomath}\begin{equation}\label{varphiFL<varphi}
	\varphiFL(x;\pi_\varphi,\rh_\varphi) \le \varphi(x) \;\text{ for every }\; x\in[0,\Pinfphi]\,,
\end{equation}\end{linenomath}
then the probability of extinction by generation $n$ satisfies
\begin{linenomath}\begin{equation}\label{PnFLgam}
	 \frac{\Pinfphi(1-\gaphi^n)}{1-\gaphi^n\Pinfphi} \le \Pnphi \le \Pinfphi \;.
\end{equation}\end{linenomath}
Equivalently, the probability $\Snphi$ of survival up to generation $n$ satisfies
\begin{linenomath}\begin{equation}\label{SnFLgam}
	\Sinfphi\le \Snphi\le \frac{\Sinfphi}{1 - \gaphi^n(1-\Sinfphi)}\,.  
\end{equation}\end{linenomath}
\end{proposition}

The key of applying this result to a given offspring distribution $\varphi$ is of course the establishment of \eqref{varphiFL<varphi}. 
These bounds yield the correct rate of approach to $\Pinfphi$ and $\Sinfphi$. However, in general, they yield little detailed information because typically $\Pinfphi$ and $\gaphi$ cannot be evaluated analytically (even for simple distributions, such as binomial or negative binomial). One remedy is to use accurate approximations for $\Pinfphi$ and $\gaphi$, which is possible for many families of distributions (see Section \ref{sec:Bounds for Sinf}).

\begin{remark}\label{rem:bound_Psurv_gamma}
If instead of \eqref{varphiFL<varphi},
\begin{linenomath}\begin{equation}\label{varphiFL>varphi}
	\varphiFL(x;\pi_\varphi,\rh_\varphi) \ge \varphi(x) \;\text{ for every }\; x\in[0,\Pinfphi]
\end{equation}\end{linenomath}
is satisfied, then 
\begin{linenomath}\begin{equation}\label{PnFLgam_reversed}
	\Pnphi \le \frac{\Pinfphi(1-\gaphi^n)}{1-\gaphi^n\Pinfphi} 
\end{equation}\end{linenomath}
and
\begin{linenomath}\begin{equation}\label{SnFLgam_reversed}
	\Snphi\ge \frac{\Sinfphi}{1 - \gaphi^n(1-\Sinfphi)} 
\end{equation}\end{linenomath}
hold. Again, this follows from \eqref{Senetas_inequ_assump+} and \eqref{Senetas_inequ}.
\end{remark}

By construction, these bounds provide excellent approximations for $\Pnphi$ and $\Snphi$ if $n$ is large, but not necessarily if $n$ is small because $\varphiFL(0)$ may differ considerably from $\varphi(0)$. Relative errors are displayed in Fig.~\ref{fig_SnGP} for a versatile class of generalized Poisson distributions.

From \eqref{SnFLgam}, we can derive a simple bound for the minimum time $T_\varphi(\ep)$ such that
\begin{linenomath}\begin{equation}\label{Tep_gen_def}
	\Snphi \le (1+\ep)\Sinfphi \; \text{for every } n \ge T_\varphi(\ep) \,. 
\end{equation}\end{linenomath}
Indeed, from \eqref{Tep_FL}, \eqref{gaFL=1/mFL}, and \eqref{SnFLgam}, we obtain $T_\varphi(\ep)\le T_{\rm FL}(\ep)$. By the construction of $\varphiFL$ in Proposition \ref{varphiFL<varphi} this yields
\begin{linenomath}\begin{equation}\label{Tep_gen_ineq}
	T_\varphi(\ep) \le \frac{\ln\Bigl(\bigl(1+\frac{1}{\ep} \bigr) \Pinfphi \Bigr)}{-\ln\gaphi}  \,.
\end{equation}\end{linenomath}
For sufficiently small $\ep$, $T_\varphi(\ep)$ is the time after which extinction of the mutant can be ignored.
We study simple approximations as well as their accuracy in Sect.~\ref{subsect:Tep}. This estimate of $T(\ep)$ will play a key role in Sect.~\ref{subsec:spread_favorable_mutant}.

In Section \ref{sec:Specific bounds} we investigate the validity of \eqref{varphiFL<varphi} for well-known families of offspring distributions. In some cases, \eqref{varphiFL<varphi} is valid for every  $x\in[0,1]$. A necessary condition for this is 
\begin{linenomath}\begin{equation}\label{mphi*gaphi<1}
	\mphi\gaphi < 1\,.
\end{equation}\end{linenomath}
Indeed, by our construction of $\varphiFL$ and by \eqref{gaFL=1/mFL}, i.e., because $\mphi^{-1} = \gaFL=\gaphi$, we obtain that \eqref{mphi*gaphi<1} holds if and only if $\varphiFL'(1)>\varphi'(1)$. The latter implies $\varphiFL(x)<\varphi(x)$ for $x$ slightly smaller than 1.

\subsection{Alternative approaches}\label{sec:alternatives}
As noted by a reviewer, a simple general upper bound for $\Pnphi$ is obtained by using concavity of the pgf $\varphi$ and starting with the observation $\Pinfphi-\Pnphi = \varphi(\Pinfphi)-\varphi(P^{n-1}_\varphi) \le \gaphi(\Pinfphi-P^{n-1}_\varphi)$. Then iteration yields
\begin{equation}\label{Simp_rate_convergence}
	\frac{\Pinfphi-\Pnphi}{\gaphi^n}  \le \Pinfphi\,, \; n\ge1 \,.
\end{equation}
Our bound \eqref{PnFLgam} yields 
\begin{equation}\label{RB_rate_convergence}
	\frac{\Pinfphi-\Pnphi}{\gaphi^n} \le \frac{\Pinfphi(1-\Pinfphi)}{1-\gaphi^n\Pinfphi}\,,
\end{equation}
where the right side converges to $\Pinfphi(1-\Pinfphi)$ as $n\to\infty$.
This yields a much tighter upper bound than \eqref{Simp_rate_convergence}, especially in the slightly supercritical case when $1-\Pinfphi=O(s)$ if $m=1+s$. Also \eqref{Simp_rate_convergence} entails a much higher estimate for $T_\varphi(\ep)$ than \eqref{Tep_gen_ineq}, which is important for the applications in Sect.~\ref{subsec:spread_favorable_mutant}. For the Poisson distribution, relative errors of these bounds and those discussed below are given in Table \ref{table_rel_err_Pollak_Agresti_etc}.

\subsubsection*{\citetalias{Pollak1971} bounds} 
For generating functions $\varphi$ with $m>1$ and $\Pinfphi>0$, \citet[][pp.~16,17]{Harris1963} proved that a constant $d>0$ exists such that, in our notation,
\begin{equation}\label{Harris_asymptotics}
	\Pnphi = \Pinfphi - d \gaphi^n + O(\gaphi^{2n})\,,
\end{equation}
\citet{Pollak1971} derived a method to obtain upper and lower bounds for $d$. His method is based on a recursive formula for $(\Pinfphi-\Pnphi)^{-1}$ that invokes series expansion of $\varphi$ about $\Pinfphi$ up to second and third order for the upper and lower bound, respectively. Application of his method requires the verification of two complicated inequalities (one for each bound) on the generating function (his two-sided inequality (2.2)). He verified both inequalities for Poisson distributions satisfying $m\PinfPoi<2$ (in fact, $m\PinfPoi=\gaPoi\le1$ holds always by Lemma \ref{lem:gaPoi_ineq}), and for negative binomial distributions satisfying two conditions. (With the help of the expansions in Sect.~\ref{subsec:NB_bounds}, it is readily shown that they are fulfilled if $m\ge1$.) \citet[][eq.~(5.2) with $r=0$]{Pollak1971} proved that 
\begin{subequations}\label{Pollak_rate_convergence}
\begin{equation}\label{Pollak_rate_convergence_a}
	\frac{\Pinfphi-\Pnphi}{\gaphi^n} \le \bar d^{(n)} \,,
\end{equation}
where
\begin{equation}\label{Pollak_bar_d^n}
	\bar d^{(n)} = \frac{2(1-\gaphi)\Pinfphi}{2(1-\gaphi) + \varphi''(\Pinfphi)\Pinfphi(1-\gaphi^n)/\gaphi}\,.
\end{equation}
\end{subequations}
At least for the Poisson distribution, Pollak's $\bar d^{(n)}$ in \eqref{Pollak_bar_d^n} is slightly smaller than our bound in \eqref{RB_rate_convergence}. In Sect.~\ref{sec:Poisson_comparison} we compare these bounds in more detail, especially numerically. The corresponding lower bound for $(\Pinfphi-\Pnphi)/\gaphi^n$ is more complicated and invokes $\varphi'''(\Pinfphi)$. We note that Pollak also derived bounds for the subcritical case.

\subsubsection*{\citetalias{Agresti1974} bounds}
\citet{Seneta1967} applied his method (Sect.~\ref{Seneta's method}) to obtain bounds for the generating function of the Poisson distribution in the subcritical case. \citet{Agresti1974} refined this approach considerably. For the subcritical case, he derived best possible fractional linear lower and upper bounds for generating functions of the form $p_0+p_1x+p_2x^2$ and $p_0+p_kx^k$ for some $k\ge1$. Agresti used these to obtain fractional linear bounds for rather general generating functions. In general, those are not best possible. By a different procedure, he derived best possible fractional linear bounds for Poisson generating functions with mean $m<1$. 

Agresti noted that the dual relation $\varphi_{\rm sub}(x)=\varphi(\Pinfphi x)/\Pinfphi$ can be used to derive bounds for a supercritical generating function $\varphi$ from bounds for the subcritical case.  
Application of Agresti's results to a supercritical $\varphi$ (with $0<\Pinfphi<1$) requires to first determine the bounds $\varphi_{\rm sub,L}$ and $\varphi_{\rm sub,U}$ for the dual subcritical $\varphi_{\rm sub}(x)$, i.e.,
\begin{equation}
	\varphi_{\rm sub,L}(x) \le \varphi_{\rm sub}(x) \le\varphi_{\rm sub,U}(x)\,.
\end{equation}
This yields the following bounds for the given $\varphi(x)$:
\begin{equation}\label{Agresti_supercrit_bounds}
	\Pinfphi\varphi_{\rm sub,L}(x/\Pinfphi) \le \varphi(x) \le \Pinfphi\varphi_{\rm sub,U}(x/\Pinfphi) \;\text{ if }\; 0\le x \le \Pinfphi\,.
\end{equation}
Clearly, equality holds in \eqref{Agresti_supercrit_bounds} if $x=\Pinfphi$. At $x=\Pinfphi$, also the first derivatives of the three functions coincide, as do the second derivatives of the first two functions. The reason is that his method for the subcritical case requires that the bounding fractional linear pgfs have the same first derivatives at $x=1$ as the given pgf $\varphi$. For the lower bound, on which we concentrate, also the second derivatives at $x=1$ must coincide. His upper bound instead satisfies the requirement that its value at $0$ equals that of the given (subcritical) pgf.
Therefore, his lower bound shares the property of yielding the correct rate of approach of $\Pnphi\to\Pinfphi$ and of $\Snphi\to\Sinfphi$ with our lower bound and with that of Pollak.  We note that the bounds in \eqref{Agresti_supercrit_bounds} are not generating functions because at $x=1$ they exceed 1.

Agresti's method has a disadvantage that may be prohibitive for applications to other generating functions. In addition to fitting the derivative at $\Pinf$ (as in our approach), it requires the determination of the supremum and infimum (with respect of $x$) of
\begin{equation}
	v(x,m) = \frac{\varphi(\Pinfphi x)/\Pinfphi-1+\gaphi(1-x)}{x \varphi(\Pinfphi x)/\Pinfphi-x+\gaphi(1-x)}\,,
\end{equation}
where this is already transformed from his version using duality, so that here $\varphi$ has mean $m>1$ and his $\la=\varphi_{\rm sub}'(1)=\gaphi$.
The parameter $\pi$ of the bounding fractional linear function (not a pgf!) is $\pi=\sup_{0\le x<1}v(x,m)$ for the upper bound and $\pi=\inf_{0\le x<1}v(x,m)$ for the lower bound. These are readily determined for $\varphiPoi$ because then Agresti proved that in this case $v(x,m)$ is strictly monotone decreasing if $x\in[0,1)$. For other generating functions, this may be much more difficult or impossible to establish. For instance, for the generalized Poisson distribution treated in Sect.~\ref{subsec:GP}, the resulting function $v(x;m,\la)$ can be decreasing in $x$ (for sufficiently small $\la$), increasing (for sufficiently large $\la$), or have a local minimum at some $x\in (0,1)$ for a small range of intermediate values $\la$.

For the Poisson distribution with $m<1$, Agresti derived the bounds $\varphi_{\rm sub,L}$ and $\varphi_{\rm sub,U}$ explicitly. We treat the resulting lower bounds for $\PnPoi$ in the supercritical case in Sect.~\ref{sec:Poisson_comparison}, where we also compare the accuracy of the bounds discussed here.

Finally we note that \citet{SagitovLindo2016} introduced a class of so-called power-fractional generating functions. Similar to linear fractional generating functions, this class has the property that it is invariant under iterations. This class is parameterized by four parameters, thus much more flexible than the linear fractional class. Branching processes with power-fractional offspring distributions were recently studied by \citet{AlsmeyerHoang2025}. It would be of interest to investigate if this class can be used to derive either more accurate bounds than the fractional linear class, or accurate bounds for distributions where the fractional linear class does not provide bounds (such as in a parameter region for distributions with at most three offspring; cf.~Sect.~\ref{sec:F3}).

\section{Bounds for the survival probabilities $\Snphi$ for common families of offspring distributions}\label{sec:Specific bounds}
In Sections \ref{sec:Poisson}, \ref{sec:Bin}, and \ref{sec:NB}, we prove validity of \eqref{varphiFL<varphi} for the families of Poisson, binomial, and negative binomial distributions, respectively. Consequently, the lower bound \eqref{PnFLgam} for $\Pnphi$ and the upper bound \eqref{SnFLgam} for $\Snphi$ are established for these distributions.  
In Section \ref{sec:F3}, we study distributions with at most three offspring and characterize when either \eqref{varphiFL<varphi} or its converse or none of both holds. Proofs are relegated to the appendix, except for the Poisson distribution for which the proof is simple enough so that the basic ideas are not hidden behind technical details. 

\subsection{The Poisson distribution}\label{sec:Poisson}
The main goal here is to prove that \eqref{varphiFL<varphi} holds for the Poisson distribution. Indeed, we will show that the inequality holds for every $x\in[0,1]$. We start by recalling some facts about the Poisson distribution with mean $m>1$. Its pgf is
\begin{linenomath}\begin{equation}
	 \varphiPoi(x;m) = e^{-m(1-x)}\,.
\end{equation}\end{linenomath}
We will need the Lambert function, or the product logarithm, $W(z)=\ProdLog(z)$, which is defined as the principal branch of the solution $w$ of $z=we^w$, $z\ge-e^{-1}$. (Lambert's $W$ function is treated in considerable detail in \citealt{Corless_etal1996} and \citealt{Wikipedia_Lambert_W}.) We will need $W(z)$ for values $z\in[-e^{-1},0]$, for which it is monotone increasing from $-1$ to 0 and concave. Then the extinction probability is
\begin{linenomath}\begin{equation}\label{PinfPoi}
	\PinfPoi(m)=-\frac{1}{m} W(-m e^{-m})\,,
\end{equation}\end{linenomath}
and 
\begin{linenomath}\begin{equation}\label{gaPoi}
	\gaPoi(m) = \varphiPoi'(\PinfPoi(m)) = m \PinfPoi(m) = - W(-m e^{-m}) \,.
\end{equation}\end{linenomath}

To establish \eqref{varphiFL<varphi}, we proceed as in the derivation of Proposition \ref{prop:bound_Psurv_gamma} and choose the parameters $\pi$ and $\rh$ of our candidate for a bounding fractional linear pgf according to \eqref{pr_solution_gam} with $a_1=\PinfPoi(m)$ and $a_2=\gaPoi(m)$. We denote them by $\pi_m$ and $\rh_m$ to indicate their dependence on $m$. Straightforward algebra yields 
\begin{linenomath}\begin{equation}\label{subst_pmrm}
	\rh_m = -\frac{W(-m e^{-m})^2+W(-m e^{-m})}{m-W(-m e^{-m})^2} \;\text{ and }\; \pi_m =-\frac{m \rh_m}{W(-m e^{-m})}\,,
\end{equation}\end{linenomath}
where the latter follows from \eqref{PinfFL} and \eqref{PinfPoi}.

\begin{theorem}\label{thm:Poi_ge_FL}
For every $m>1$, the pgfs $\varphiPoi(x;m)$ and $\varphiFL(x;\pi_m,\rh_m)$ satisfy
\begin{linenomath}\begin{equation}\label{Poi_ge_FL}
	\varphiFL(x;\pi_m,\rh_m) \le \varphiPoi(x;m) \;\text{ for every } x\in[0,1]\,.
\end{equation}\end{linenomath}
Equality holds only at $x=\PinfPoi$ and $x=1$.
\end{theorem}

Proposition \ref{prop:bound_Psurv_gamma} immediately yields
\begin{corollary}\label{cor:Poi_ge_FL}
Given a Poisson offspring distribution with mean $m>1$, the probability of extinction by generation $n$, $\PnPoi$, satisfies the inequality \eqref{PnFLgam}, and the probability of survival up to generation $n$, $\SnPoi$, satisfies the inequality \eqref{SnFLgam}, each with $\varphi=\varphiPoi$.
\end{corollary}

In the proof of Theorem \ref{thm:Poi_ge_FL} we will need some inequalities.

\begin{lemma}\label{lem:gaPoi_ineq}
The following inequalities hold:
\begin{linenomath}\begin{equation}\label{gaPoi_ineq}
	2- m < \gaPoi(m)< \frac{1}{m}
\end{equation}\end{linenomath}
and
\begin{linenomath}\begin{equation}\label{PinfPoi_ineq}
	\frac{2}{m}-1 < \PinfPoi(m)< \frac{1}{m^2} \,,
\end{equation}\end{linenomath}
where both lower bounds are 0 if $m\ge2$.
\end{lemma}

\begin{proof}
Because of \eqref{gaPoi} it is sufficient to prove \eqref{gaPoi_ineq}.
We start with the right hand side. By \eqref{gaPoi} and the definition of $W(z)$, $\gaPoi$ satisfies
\begin{linenomath}\begin{equation}\label{gaPoi_relation}
	-m e^{-m} = -\gaPoi(m) e^{-\gaPoi(m)}\,.
\end{equation}\end{linenomath}
By the properties of $W$ we have $\gaPoi(1)=1$, $0<\gaPoi(m)<1$ if $m>1$, and $\gaPoi(m)$ decreases monotonically to 0 as $m\to\infty$.
The function $g(x)=-x e^{-x}$ is monotone decreasing from 0 at $x=0$ to $-e^{-1}$ at $x=1$. Therefore, $\gaPoi < \frac{1}{m}$ if and only if $g(\gaPoi)> g(1/m)$, which is equivalent to $-m e^{-m} > -\frac{1}{m} e^{-1/m}$ by using \eqref{gaPoi_relation}. The latter inequality can be rearranged to $m^2 e^{1/m-m} < 1$, which is easily verified because the left hand side equals 1 if $m=1$ and is monotone decreasing in $m$.

To prove the left hand side of \eqref{gaPoi_ineq}, we show $W(-me^{-m})<m-2$.
If $m\ge2$, this is trivially satisfied because $W(-me^{-m})<0$ whenever $m\ge1$. If $m<2$
we use that $0<x<\Pinf$ is equivalent to $\varphi(x)> x$ for any generating function. With $x=\frac{2}{m}-1$, this shows that $\PinfPoi= -m^{-1} W(-m e^{-m})> \frac{2}{m}-1$ if and only if $\varphiPoi(\frac{2}{m}-1) = e^{2-2m} > \frac{2}{m}-1$. The latter inequality is readily confirmed, for instance by showing that $\frac{d}{dm}\bigl( e^{2-2m}/(\frac{2}{m}-1)\bigr) = 2e^{2-2m}(m-1)^2/(m-2)^2 >0$.
\end{proof}

\begin{proof}[Proof of Theorem \ref{thm:Poi_ge_FL}]
Proving \eqref{Poi_ge_FL} is equivalent to showing 
\begin{linenomath}\begin{equation}\label{fP}
	\fP(x) = \ln\varphiPoi(x;m) - \ln\varphiFL(x;\pi_m,\rh_m) \ge 0 \;\text{ for every  }\; x\in[0,1]\,,
\end{equation}\end{linenomath}
where we omit the dependence of $\fP$ on $m$.
Proving \eqref{fP} is simplified by the fact that $\ln\varphiPoi(x;m)=m(x-1)$. We easily infer from the properties of $\varphiPoi$ stated above and the definition of $\varphiFL(x;\pi_m,\rh_m)$ that $\fP(\PinfPoi) = 0$, $\fP(1)=0$, $\fP'(\PinfPoi)=0$, and $\fP'(1) = m-\gaPoi^{-1}<0$, where we used \eqref{gaFL=1/mFL} and the right-hand side of \eqref{gaPoi_ineq}. A typical graph of $\fP$ is shown in Fig.~\ref{fig_fPoi}.  

Now we show $\fP''(\PinfPoi) = -(\ln\varphiFL)''(\PinfPoi) > 0$. For a general fractional linear pgf we get $-(\ln\varphiFL)''(\PinfFL) = \frac{(1-\pi)\pi^2(1-\pi-2\rh)}{(1-\rh)^2\rh^2}$ because $\PinfFL=\rh/\pi$. This is positive if and only if $\pi+2\rh<1$. With $\pi =\pi_m$ and $\rh =\rh_m$ from \eqref{subst_pmrm}, we obtain after some calculation that $\pi_m+2\rh_m<1$ if and only if 
\begin{linenomath}\begin{equation}\label{ineq_fP''}
	\frac{\bigl(m-2W(-m e^{-m})\bigr) \bigl(1+W(-m e^{-m})\bigr)}{m-W(-m e^{-m})^2} < 1\,.
\end{equation}\end{linenomath}
Each of the three factors on the left hand side is positive if $m>1$ because $-1< W(-m e^{-m})<0$. Therefore rearrangement shows that the inequality \eqref{ineq_fP''} is equivalent to $W(-m e^{-m}) < m - 2$, which we proved in Lemma \ref{lem:gaPoi_ineq}.

\begin{figure}[t!]
\centering
\includegraphics[width=0.6\textwidth]{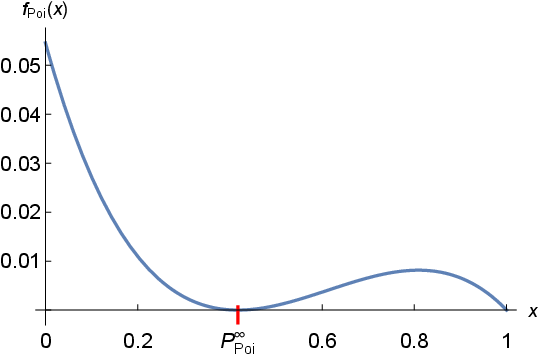}
\caption{The graph of the function $\fP(x)$ with $m=1.5$. Then $\pi_m\approx0.506$, $\rh_m\approx0.211$, $\PinfPoi\approx0.4172$. As $m$ decreases to 1, $\PinfPoi$ increases to 1, and $(\pi_m,\rh_m)$ approaches $(\tfrac13,\tfrac13)$.}
\label{fig_fPoi}
\end{figure}

If we can show that $\fP'''(x) = -(\ln\varphiFL)'''(x)<0$ for every $0<x<1$, then $\fP'(x)$ is strictly concave and has exactly one zero between $\PinfPoi$ and 1, because $\fP'(\PinfPoi)=0$, $\fP''(\PinfPoi)>0$, and  $\fP'(1)<0$. Because $\fP(\PinfPoi)=\fP(1)=0$ and $\fP'(1)<0$, $\fP$ must have a local maximum at this zero. It also follows that $\fP'(x)<0$ if $0<x<\PinfPoi$. Hence, $\fP(x)\ge 0$ for every $0\le x\le 1$, with equality only at $x=\PinfPoi$ and $x=1$.

It remains to determine the sign of $\fP'''(x)$. It is straightforward to check (Sect.~2.2 in the supplementary {\it Mathematica} notebook) that
\begin{linenomath}\begin{equation}
	(\ln\varphiFL)'''(x;\pi,\rh) = \dfrac{2(1-\pi)(1-\rh)d(x)}{\bigl(\rh(1-x)+(1-\pi)x\bigr)^3(1-\pi x)^3}\,,
\end{equation}\end{linenomath}
where
\begin{linenomath}\begin{equation}
	d(x) = \frac{3}{4}\bigl(2x\pi(1-\pi-\rh) - (1-\pi-\rh-\rh\pi)\bigr)^2 + \frac{1}{4}(1-\pi)^2(1-\rh)^2\,.
\end{equation}\end{linenomath}
Therefore, $(\ln\varphiFL)'''(x;\pi,\rh)>0$ for every $x\in(0,1)$ and every admissible pair $(\pi,\rh)$. This finishes the proof.
\end{proof}

We note that the proof implies that $e^{-m}=\varphiPoi(0;m) > \varphiFL(0;\pi_m,\rh_m)=\rh_m$, which is not easy to establish directly. It also implies that $\varphiPoi'(0)<\varphiFL'(0)$.

In the proof above we showed that $\fP'(x)<0$ if $x\in[0,\PinfPoi)$. By a simple calculation we infer that the maximum of the relative error $(\varphiPoi(x)-\varphiFL(x))/\varphiPoi(x)$ on the interval $[0,\PinfPoi]$ is achieved at $x=0$.
At $x=0$, we find $\lim_{m\downarrow1}\varphiFL(0;\pi_m,\rh_m)= \tfrac13 < \tfrac1e = \lim_{m\downarrow1}\varphiPoi(0;m)$, which yields a relative error of $\approx0.094$ if $m=1$. 
If $m=1+s$ and $s$ is small, such as $s<0.3$, an accurate approximation of the maximum relative error is $1-\tfrac{e}3-\tfrac{e}{27}s \approx 0.094-0.101s$. The relative error decreases to 0 as $m\to\infty$ (unsurprisingly, because $\PinfPoi\to0$).
Therefore, the upper bound in \eqref{SnFLgam} is not an accurate approximation if $m$ is close to 1 and $n$ is small (see the case $\la=0$ in Fig.~\ref{fig_SnGP}). 

Series expansions in $s$ of $\SinfPoi$ and $\gaPoi$ when $m=1+s$, as well as upper and lower bounds of $\SinfPoi$, are presented in Sect.~\ref{subsec:Poisson_bounds}.

\subsubsection{Comparison of Pollak's and Agresti's bounds with our and the simple bound}\label{sec:Poisson_comparison}
For the Poisson distribution, the term $ \varphi''(\Pinfphi)\Pinfphi(1-\gaphi^n)/\gaphi$ in Pollak's bound $\bar d^{(n)}$ in \eqref {Pollak_bar_d^n} for $(\Pinfphi-\Pnphi)/\gaphi^n$  simplifies to $\gaPoi(1-\gaPoi^n)$. Therefore,
\begin{equation}\label{Pollak_Agresti_bound}
	\bar d^{(n)} = \frac{\PinfPoi}{1+\frac{\gaPoi(1-\gaPoi^n)}{2(1-\gaPoi)}}\,,
\end{equation}
whence
\begin{equation}\label{bar d}
	\bar d = \lim_{n\to\infty}\bar d^{(n)} = \frac{\PinfPoi}{1+\tfrac12\gaPoi(1-\gaPoi)^{-1}}
\end{equation}
is an upper bound for Harris' constant $d$ in \eqref{Harris_asymptotics}. If $m=1+s$, then $\gaPoi = 1-s + \tfrac23s^2 + O(s^3)$ (Sect.~\ref{subsec:Poisson_bounds}), so that $\left(1+\tfrac12\gaPoi(1-\gaPoi)^{-1}\right)^{-1} = 2s - \tfrac{10}{3}s^2 + O(s^3)$. We recall from \eqref{RB_rate_convergence} that our corresponding upper bound for $d$ is $\PinfPoi(1-\PinfPoi)$, where $1-\PinfPoi= 2s -\tfrac83s^2 + O(s^3)$. Thus, Pollak's bound is slightly more accurate.

As already noted in Sect.~\ref{sec:alternatives}, \citet{Agresti1974} derived upper and lower bounds for the Poisson distribution with $m<1$. From the parameters of his bounds for the subcritical case, the parameters $\pi$ and $\rh$ for the  fractional linear bounds in the supercritical case can be computed and are given in Sect.~2.3 of the supplementary {\it Mathematica} notebook. It turns out that the lower bounds for $\PnPoi$ of Pollak and Agresti coincide (as already noted by Agresti). Agresti states that his upper bound, which is also given in the supplementary notebook, performs favorably compared with Pollak's. 

Table \ref{table_rel_err_Pollak_Agresti_etc} shows the relative errors $(\Sn_{\rm app}-\SnPoi)/\SnPoi$ produced by the bounds discussed above, where $\SnPoi$ is the exact value obtained by iteration of the pgf $\varphiPoi$. The data confirm that Pollak's and Agresti's bounds perform slightly better than ours. The reason is that both are based on fitting also the second derivative at $\PinfPoi$, i.e., Agresti's bounding fractional linear function has the property that it not only coincides with $\varphiPoi$ at $\PinfPoi$, but also its first and second derivative do. Therefore, it cannnot be a pgf in the supercritical case (indeed it exceeds 1 at $x=1$). Our method posits a bounding fractional linear pgf, whence only it and its first derivative can be fitted at $\PinfPoi$. Pollak does not construct bounding functions.

\begin{table}[t!]
\centering
\begin{tabular}{l|ccccccc}
\hline\hline
	& 		$n=1$ & $n=5$ & $n=10$ & $n=20$ & $n=50$ & $n=100$  \\
	\hline 
	$m=1.5$  \\
	\eqref{Simp_rate_convergence}
		& 0.0863 & 0.0295 & 0.00307 & $2.8\times10^{-5}$ & $2.2\times10^{-11}$ & 0  \\
	\eqref{SnFLgam} 
		& 0.0153 & 0.0035 & 0.00034 & $3.1\times10^{-6}$ & $2.5\times10^{-12}$ & 0  \\
	\eqref{Pollak_Agresti_bound}
		& 0.0062 & 0.0010 & 0.00009 & $8.3\times10^{-7}$ & $6.4\times10^{-13}$ & 0 \\
\hline 
	$m=1.1$  \\
	\eqref{Simp_rate_convergence}
		& 0.3832 & 0.9869 & 0.94154 & 0.47327 & 0.02823 & $2.1\times10^{-4}$  \\
	\eqref{SnFLgam} 
		& 0.0420 & 0.0372 & 0.02084 & 0.00705 & 0.00035 & $2.5\times10^{-6}$  \\
	\eqref{Pollak_Agresti_bound}
		& 0.0342 & 0.0262 & 0.01321 & 0.00404 & 0.00019 & $1.4\times10^{-6}$  \\
\hline 
	$m=1.02$  \\
	\eqref{Simp_rate_convergence}
		& 0.5343 & 2.2140 & 3.7106 & 5.405 & 5.589 & 2.802  \\
	\eqref{SnFLgam} 
		& 0.0518 & 0.0587 & 0.04439 & 0.02777 & 0.01050 & 0.00316 \\
	\eqref{Pollak_Agresti_bound}
		& 0.0497 & 0.0546 & 0.03994 & 0.02386 & 0.00827 & 0.00233 	\\
	\hline\hline
\end{tabular}
\caption{The table shows the relative errors $(\Sn_{\rm app}-\SnPoi)/\SnPoi$ of the upper bound $\Sn_{\rm app}$ for $\SnPoi$ obtained from the simple method in \eqref{Simp_rate_convergence}, from \eqref{SnFLgam}, and from Pollak's   \eqref{Pollak_Agresti_bound} (and Agresti's equivalent bound). The relative error of the simple bound tends to 0 very slowly. If $m=1.02$, it takes 506 generations to decrease below 0.001.}\label{table_rel_err_Pollak_Agresti_etc}
\end{table}

\FloatBarrier

\subsection{The binomial distribution}\label{sec:Bin}
The binomial distribution has the pgf 
\begin{linenomath}\begin{equation}
	\varphiBin(x;n,p) = (1-p+px)^n\,.
\end{equation}\end{linenomath}
We assume $n\ge2$ and $\mBin=np>1$. 
Let $\PinfBin$ denote the extinction probability, i.e., the unique solution of $\varphiBin(x)=x$ in $(0,1)$. We set
\begin{linenomath}\begin{equation}
	\xi = ({\PinfBin})^{1/n}
\end{equation}\end{linenomath}
and $p=\dfrac{1-\xi}{1-\xi^n}$, and parameterize $\varphiBin$ by $n$ and $\xi$.
Then 
\begin{equation}\label{phiBin_xi}
	\varphiBin\Bigl(x;n,\frac{1-\xi}{1-\xi^n}\Bigr) = \left(\frac{x(1-\xi)+(\xi-\xi^n)}{1-\xi^n}\right)^n\,.
\end{equation}

From eqs.~\eqref{pr_ab1} and \eqref{pr_solution_gam}, we infer that the fractional linear pgf $\varphiFL(x;\pBin,\rBin)$ with the parameters
\begin{equation}
	\pBin = \frac{\xi(1-\xi^n + n \xi^n)-n \xi^n}{\xi(1-\xi^n + n \xi^{2n})-n \xi^{2n}}\;\text{ and } \rBin = \xi^n \,\pBin
\end{equation}
has the same extinction probability, $\PinfBin$, and the same rate of convergence,
\begin{equation}\label{gaBin_xi}
	\gaBin = \frac{n \xi^n(1-\xi)}{\xi(1-\xi^n)} \,,
\end{equation}
as the binomial.

In Appendix \ref{sec:Proof_Bin}, we prove

\begin{theorem}\label{thm:Bin_ge_FL}
For every $n\ge2$ and every $\xi\in(0,1)$, the pgfs $\varphiBin$ and $\varphiFL$ satisfy
\begin{linenomath}\begin{equation}\label{Bin_ge_FL}
	\varphiFL(x;\pBin,\rBin) \le \varphiBin\Bigl(x;n,\frac{1-\xi}{1-\xi^n}\Bigr) \;\text{ for every } x\in[0,1]\,.
\end{equation}\end{linenomath}
Equality holds if and only if $x = \PinfBin$ or $x=1$.
\end{theorem}

Proposition \ref{prop:bound_Psurv_gamma} immediately yields
\begin{corollary}
Given a binomial offspring distribution with mean $\mBin>1$, the probability of extinction by generation $\tau$, $P_{\rm Bin}^{(\tau)}$, satisfies the inequality \eqref{PnFLgam}, and the probability of survival up to generation $\tau$, $S_{\rm Bin}^{(\tau)}$, satisfies the inequality \eqref{SnFLgam}, each with $\varphi=\varphiBin$.
\end{corollary}

Series expansions in $s$ of $\SinfBin$ and $\gaBin$ when $\mBin=1+s$, as well as upper and lower bounds of $\SinfBin$, are presented in Sect.~\ref{subsec:Bin_bounds}.

\subsection{The negative binomial distribution}\label{sec:NB}
The negative binomial distribution has the pgf 
\begin{linenomath}\begin{equation}
	\varphiNB(x;r,p) = \frac{p^r}{(1-(1-p)x)^r}\,.
\end{equation}\end{linenomath}
The mean and the variance are $\mNB=r\frac{1-p}{p}$, $\si^2_{\rm NB}= r\frac{1-p}{p^2}$, respectively. Because for $r=1$ a geometric distribution is obtained, and nothing remains to be proved, we assume $r\ge2$ and $\mNB>1$.
Let $\PinfNB$ denote the  extinction probability, i.e., the unique solution of $\varphiNB(x)=x$ in $(0,1)$. We set
\begin{linenomath}\begin{equation}
	\ze = ({\PinfNB})^{1/r}
\end{equation}\end{linenomath}
and $p=\dfrac{\ze(1-\ze^r)}{1-\ze^{r+1}}$, and parameterize the negative binomial distribution by $r$ and $\ze$. By our general assumptions we have $0<\ze<1$.
Then 
\begin{equation}\label{phiNB_ze}
	\varphiNB\Bigl(x;r,\frac{\ze(1-\ze^r)}{1-\ze^{r+1}}\Bigr) = \left(\frac{\ze(1-\ze^r)}{1-\ze^{r+1}-x(1-\ze)}\right)^r
\end{equation}
and $\mNB=\dfrac{r(1-\ze)}{\ze(1-\ze^r)}$. Straightforward algebra shows that
\begin{linenomath}\begin{equation}\label{gaNBzet}
	\gaNB = \frac{r (1-p) \PinfNB}{1-(1-p)\PinfNB} = \frac{r(1-\ze)\ze^r}{1-\ze^r}\,.
\end{equation}\end{linenomath}

From eqs.~\eqref{pr_ab1} and \eqref{pr_solution_gam}, we infer that the fractional linear pgf $\varphiFL(x;\pNB,\rNB)$ that has the same extinction probability, $\PinfNB$, and the same rate of convergence, $\gaNB$, as the negative binomial has the parameters
\begin{linenomath}\begin{equation}\label{pNB_rNB}
	\pNB = \frac{1-[1+r(1-\ze)]\ze^r}{1-[1+r(1-\ze)\ze^r]\ze^r}\;\text{ and } \rNB = \ze^r \,\pNB\,.
\end{equation}\end{linenomath}

In Appendix \ref{sec:Proof_NB}, we prove

\begin{theorem}\label{thm:NB_ge_FL}
For every $r\ge2$ and every $\ze\in(0,1)$, the pgfs $\varphiNB$ and $\varphiFL$ satisfy
\begin{linenomath}\begin{equation}\label{varphiNB_ge_varphiFL}
	\varphiFL(x;\pNB,\rNB) \le \varphiNB\Bigl(x;r,\frac{\ze(1-\ze^r)}{1-\ze^{r+1}}\Bigr) \;\text{ for every }x\in[0,\PinfNB]\,.
\end{equation}\end{linenomath}
Equality holds if and only if $x = \PinfNB$.
\end{theorem}

Proposition \ref{prop:bound_Psurv_gamma} immediately yields
\begin{corollary}
Given a negative binomial offspring distribution with mean $\mNB>1$, the probability of extinction by generation $n$, $\PnNB$, satisfies the inequality \eqref{PnFLgam}, and the probability of survival up to generation $n$, $\SnNB$, satisfies the inequality \eqref{SnFLgam}, each with $\varphi=\varphiNB$.
\end{corollary}

We conjecture that the inequality in \eqref{varphiNB_ge_varphiFL} holds for every $x\in[0,1]$, although our proof yields it only for a smaller interval that contains $[0,\PinfNB]$.
However, we show in Appendix \ref{sec:Proof_NB} that \eqref{varphiNB_ge_varphiFL} is valid for every $x\in[0,1]$ if $r=2,\ldots, 6$.  In addition, we prove that $\mNB \gaNB <1$, which implies $\mNB \PinfNB <1$ and that \eqref{varphiNB_ge_varphiFL} holds for $x$ sufficiently close to 1; see \eqref{mphi*gaphi<1}.
Also the convergence of the negative binomial to the Poisson distribution as $r\to\infty$ (with $\mNB$ fixed) supports our conjecture.

Series expansions in $s$ of $\SinfNB$ and $\gaNB$ when $\mNB=1+s$, as well as upper and lower bounds for $\SinfNB$, are presented in Sect.~\ref{subsec:NB_bounds}.

\subsection{Offspring distributions with at most three offspring}\label{sec:F3}
\begin{figure}[th]
\vspace{-6mm}
\centering
\begin{tabular}{cc}
\includegraphics[width=0.43\textwidth]{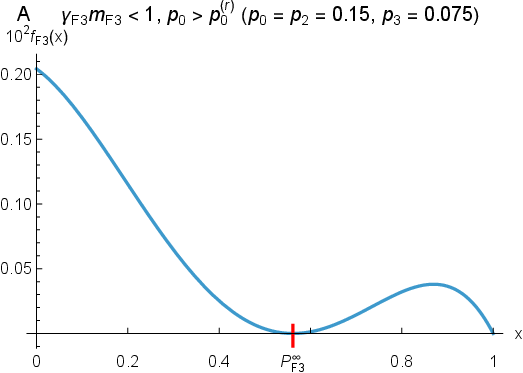} & \includegraphics[width=0.43\textwidth]{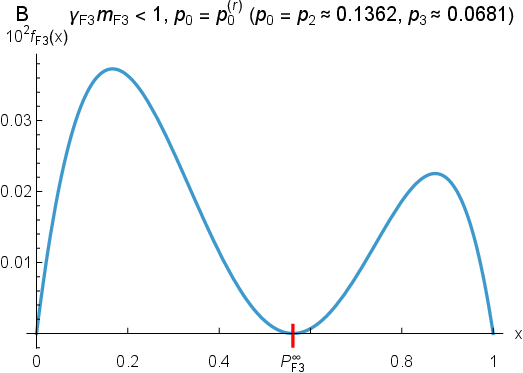} \vspace{2mm} \\ 
\includegraphics[width=0.43\textwidth]{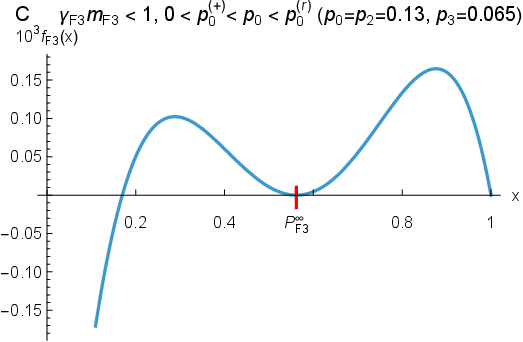} & \includegraphics[width=0.43\textwidth]{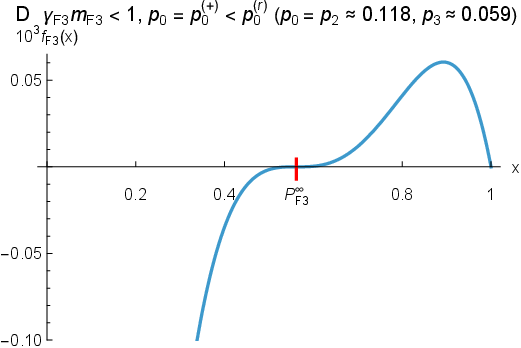} \vspace{2mm} \\ 
\includegraphics[width=0.43\textwidth]{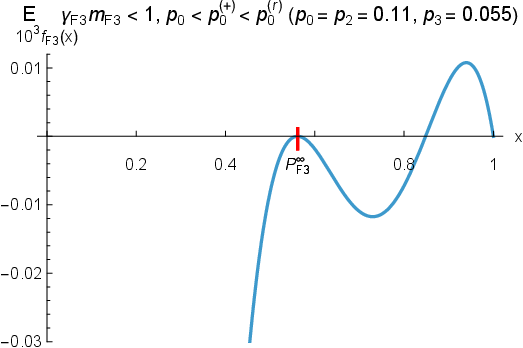} & \includegraphics[width=0.43\textwidth]{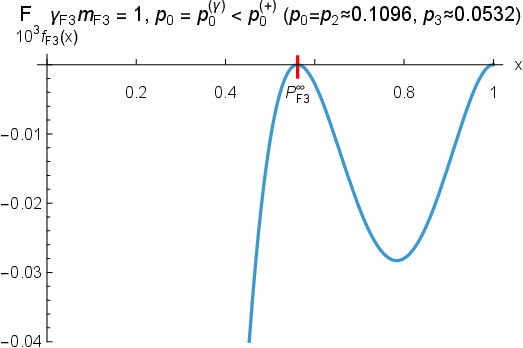} \vspace{2mm} \\
\includegraphics[width=0.43\textwidth]{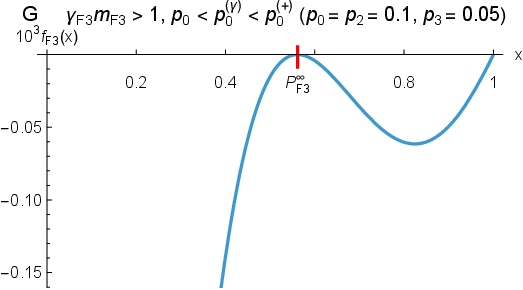}
\end{tabular}
\caption{\footnotesize{Possible shapes of graphs of $\fFth(x)=\varphiFth(x)-\varphiFL(x)$. All possible cases are obtained by choosing $p_0=p_2$, $p_3=\tfrac12p_2$, and varying $p_2$. Then the relations between $p_0$, $p_2$, and $p_3$ are retained as $p_2$ or $p_1=1-\tfrac52p_2$ varies (e.g., $p_1=0.625$ in panel A). In the degenerate case of panel D, we have $p_0=\pnplus<\pnr$, so that $\fFth''(\PinfFth)=0$ and $\fFth'''(\PinfFth)>0$. In the degenerate case of panel F, we have $\fFth'(1)=0$ and $\fFth''(1)<0$. In addition to the indicated relations, $\pnr>\pnplus>0$ holds in A and B, $\pnplus>\pnr$ in F and G, and $p_0>\pnga$ in A -- E.
In all cases, $\PinfFth=\tfrac12(\sqrt{17}-3)\approx0.56155$, and $\mFth=1+p_2$. Figure A applies if $0.4\ge p_2>\tfrac12(1-3/\sqrt{17})$, and the lower bound yields the critical case B. The critical case D occurs if $p_2=\tfrac{2}{17}$, F applies if $p_2=-\tfrac12+\tfrac{5}{34}\sqrt{17}$, and G applies for all smaller values of $p_2$. The values of $\fFth(0)$ are $\approx 0.00081$, $-0.00222$, $-0.002959$, $-0.003273$, and $-0.00376$ in panels C, D, E, F, and G, respectively. Note that the vertical scales in A and B differ from those in the other panels.}
} 
\label{fig_F3}
\end{figure}

\begin{figure}[th]
\vspace{5mm}
\centering
\begin{tabular}{ll}
{\bf A} & {\bf B}  \\
\includegraphics[width=0.45\textwidth]{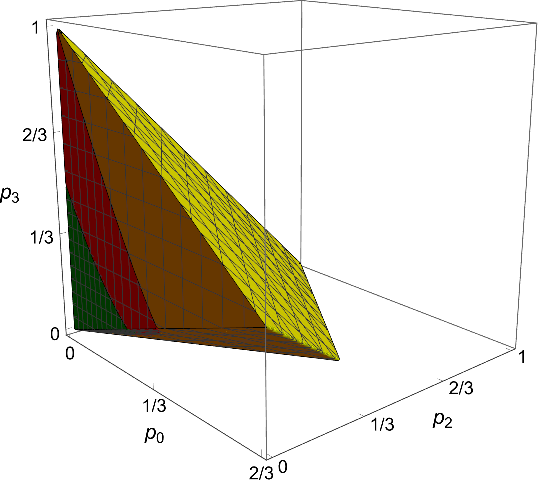} & \includegraphics[width=0.47\textwidth]{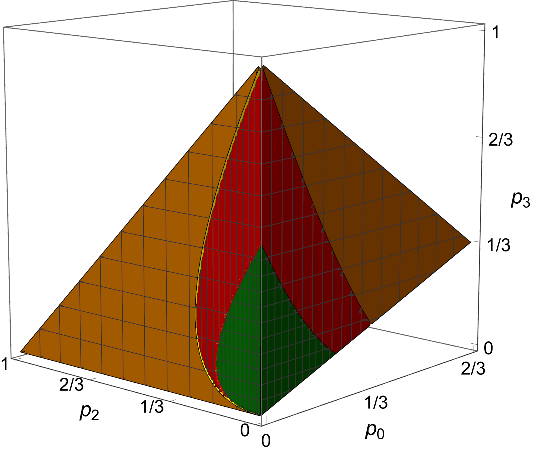} 
\end{tabular}
\caption{\footnotesize{The three regions defined in Corollary \ref{cor:main_result_F3} shown from two angles in panels A and B. The region defined by \eqref{lower_bound} is shown in shades of yellow and brown. Here, the extinction probability $\PnFth$ can be bounded from below by the fractional linear extinction probability $\PnFL$ obtained from \eqref{rFthpFth}. The yellow plane in A is the boundary $p_0+p_2+p_3=1$ ($p_1=0$). The region defined by \eqref{no_bound} is shown in shades of red. Here, $\PnFth$ cannot be bounded by $\PnFL$ from one side. The region defined by \eqref{upper_bound} is shown in shades of green. Here, $\PnFth$ is bounded from below by $\PnFL$. The boundary plane $p_0=p_2+2p_3$ ($\mFth=1$) is visible in A, close to the bottom of the cube.}
} 
\label{fig_region_F3}
\end{figure}

Here we investigate offspring distributions $\{p_k\}$ satisfying 
\begin{linenomath}\begin{equation} \label{pk_F3}
	p_0 >0, \; p_3>0,  \text{ and $p_k=0$ if $k\ge4$}\,.
\end{equation}\end{linenomath}
We exclude the trivial case $p_0=0$ and the simple case $p_3=0$, which is treated separately in Remark~\ref{rem:p3=0}. As a consequence, we assume $0<p_2+p_3<1$.
The main results are Theorem \ref{thm:main_result_F3} and Corollary \ref{cor:main_result_F3}, which provide a complete characterization when the sequence of probabilities $\Pexn$ can be bounded from below as in Proposition \ref{prop:bound_Psurv_gamma}, or from above as in Remark \ref{rem:bound_Psurv_gamma}, or neither nor. The formulation of the main results requires considerable preparation. Illustrations of the main results are shown in Figures \ref{fig_F3} and \ref{fig_region_F3}.

We express all relevant functions in terms of $p_0$, $p_2$, and $p_3$ by setting $p_1=1-p_0-p_2-p_3$. Then the pgf is
\begin{linenomath}\begin{equation} \label{varphi_F3}
	\varphiFth(x) = \varphiFth(x;p_0,p_2,p_3) = p_0 + (1-p_0-p_2-p_3)x + p_2x^2 + p_3 x^3\,,
\end{equation}\end{linenomath}
the expected number of offspring is
\begin{linenomath}\begin{equation} \label{m_F3}
	\mFth = \varphiFth'(1) = 1-p_0+p_2+2p_3\,,
\end{equation}\end{linenomath}
where $\varphiFth'$ always refers to the derivative with respect to $x$. Throughout, we assume $\mFth>1$, i.e.,  $p_0<p_2+2p_3$. The probability of (ultimate) extinction is 
\begin{linenomath}\begin{equation} \label{Pinfty_F3}
	\PinfFth = \frac{\sqrt{4p_0p_3+(p_2+p_3)^2}-(p_2+p_3)}{2p_3}\,.
\end{equation}\end{linenomath}
Our assumptions imply $0<\PinfFth<1$.  

Following \eqref{def_gamma}, we define $\gaFth =\varphiFth'(\PinfFth)$, which is the rate of convergence of $\PnFth$ to $\PinfFth$. A straightforward calculation yields 
\begin{linenomath}\begin{equation} 
	\gaFth = 1 - \frac{(p_2+3p_3)\sqrt{4p_0p_3+(p_2+p_3)^2} - 4p_0p_3 - (p_2+p_3)^2}{2p_3}\,.
\end{equation}\end{linenomath}
In the limit $p_3\to0$, we obtain $\gaFth \to 1+p_0-p_2$. With the help of Section 5 of the supplementary {\it Mathematica} notebook all formulas can be expeditiously verified.

We begin by defining the prospective bounding fractional linear pgf $\varphiFL(x;\pFth,\rFth)$. Following Proposition \ref{prop:bound_Psurv_gamma}, we require the conditions in \eqref{cond_Pinf_gainf}, i.e., $\PinfFL = \PinfFth$ and $\gaFL=\gaFth$. These hold if and only if the parameters $\pi=\pFth$ and $\rh=\rFth$ of $\varphiFL$ are 
\begin{linenomath}\begin{equation}\label{rFthpFth}
	\rFth = \frac{2p_0\sqrt{4p_0p_3+(p_2+p_3)^2}}{(p_2+p_3)+(1+2p_0)\sqrt{4p_0p_3+(p_2+p_3)^2}} \; \text{and} \;
	\pFth = \frac{\rFth}{\PinfFth} \,.
\end{equation}\end{linenomath}

Throughout, we consider the following region of admissible parameters:
\begin{linenomath}\begin{equation}\label{R}
	R=\{(p_0,p_2,p_3): p_0>0,\, p_2\ge0,\, p_3>0,\, p_0+p_2+p_3\le1,\, p_0<p_2+2p_3\} \,.
\end{equation}\end{linenomath} 
A triple $(p_0,p_2,p_3)\in R$ defines a probability distribution with $\mFth>1$ and  $0<\PinfFth < 1$. 

Our main goal will be to determine when
\begin{linenomath}\begin{equation}\label{def_f}
	\fFth(x) = \varphiFth(x;p_0,p_2,p_3) - \varphiFL(x;\pFth,\rFth) 
\end{equation}\end{linenomath}
is positive or negative. (We use properties such as positive, increasing, or convex in the strict sense.) We recall that our construction of $\varphiFL$ implies
\begin{linenomath}\begin{equation}\label{properties_f_1}
	\fFth(\PinfFth)=0 \,,\;  \fFth'(\PinfFth)=0 \,, \,\text{ and }\; \fFth(1)=0\,.
\end{equation}\end{linenomath}

We define the following quantities:
\begin{linenomath}\begin{equation}\label{pnplus}
	\pnplus = \frac{p_3-(p_2+p_3)^2}{4p_3} \,,
\end{equation}\end{linenomath}
\begin{linenomath}\begin{equation}\label{pnr}
	\pnr = \frac{1}{2}-\frac{p_2+p_3}{8p_3}\Bigl(p_2+p_3 + \sqrt{8p_3 +(p_2+p_3)^2}\Bigr) \,,
\end{equation}\end{linenomath}
and 
\begin{linenomath}\begin{equation}\label{pnga}
	\pnga = \frac{1}{2}-\frac{1}{8p_3} \Bigl(2(p_2+p_3)^2 + (p_2+3p_3)\sqrt{8p_3 +(p_2+3p_3)^2}  - (p_2+3p_3)^2 \Bigr)\,.
\end{equation}\end{linenomath}
In the following remark, the meaning of these quantities is explained.

\begin{remark}\label{rem:props_p0+_etc}
(a) $\pnplus$ is the only potentially admissible solution of $\fFth''(\PinfFth;p_0,p_2,p_3) = 0$, i.e., such that $(\pnplus,p_2,p_3) \in R$ under suitable conditions. This follows from
$\fFth''(\PinfFth) = \dfrac{1}{p_3}(4p_0p_3 + (p_2+p_3)^2-p_3)\left(p_2+3p_3-\sqrt{4p_0p_3 + (p_2+p_3)^2}\right)$ because the second factor yields $\pnplus$, and the last factor is positive if $p_0<p_2+2p_3$. We note that 
\begin{linenomath}\begin{equation}\label{pnplus_eq}
	 p_0>\pnplus \; \text{ if and only if }\; \fFth''(\PinfFth;p_0,p_2,p_3) > 0\,.
\end{equation}\end{linenomath}
Therefore, $\fFth(x)>0$ near the critical point $x=\PinfFth$ if and only if $p_0>\pnplus$.

(b) $\pnr$ is the only potentially admissible solution of $p_0 = \rFth$. We note that $\fFth(0)=p_0-\rFth$. Therefore, 
\begin{linenomath}\begin{equation}\label{pnr_eq}
	p_0>\pnr \; \text{ if and only if }\; \fFth(0)>0 \,.
\end{equation}\end{linenomath}

(c) $\pnga$ is the only potentially admissible solution $p_0$ of $\gaFth \mFth = 1$ (the proof is outlined in eq.~\eqref{gaFth*mFth-1} of the appendix; a further solution is $p_0=p_2+2p_3$, which has multiplicity two and is not admissible).  We recall from \eqref{mphi*gaphi<1} that $\gaFth \mFth \le 1 $ if and only if $\fFth'(1)\le0$.
This is a necessary condition for $\fFth(x)$ to be positive for $x\in(\PinfFth,1)$.
\end{remark}

In the following remark we summarize the admissibility conditions of and the relations between $\pnplus$, $\pnr$, and $\pnga$. 

\begin{remark}\label{rem_props p*}
(a) We have $\pnplus > 0$ if and only if  
\begin{linenomath}\begin{equation}\label{pnplus_pos}
	p_2 < \sqrt{p_3} - p_3\,.
\end{equation}\end{linenomath}
Because we assume $p_3>0$, we have $\pnplus<\tfrac14$. Note that if $\pnplus>0$ then $p_2<\tfrac14$.

(b) We have $\pnr>0$ if and only if \eqref{pnplus_pos} holds. Indeed, a simple calculation shows that $\pnr>\pnplus$ if and only if $\pnplus>0$, and $\pnr=\pnplus$ if and only if $\pnplus=0$. Straightforward algebra also shows that $\pnr+p_2+p_3<1$ always holds (because $0<p_2+p_3<1$) and that $\pnr < p_2+2p_3$ if and only if 
\begin{linenomath}\begin{equation}\label{pnr<m}
	p_3 > \frac{1}{6} \;\text{ or }\; p_2 > \tfrac12\Bigl(\sqrt{p_3(4+p_3)}-5p_3 \Bigr)\,,
\end{equation}\end{linenomath}
where the maximum value $p_2=3-2\sqrt2\approx0.172$ is assumed at $p_3=3/\sqrt2-2\approx 0.121$.

(c) We have $\pnga>0$ if and only if 
\begin{linenomath}\begin{equation}\label{pnga>0}
	p_3<\tfrac12 \;\text{ and }\; p_2< \tfrac12\left(\sqrt{p_3(4+p_3)}-3p_3\right)\,.
\end{equation}\end{linenomath}
Furthermore, 
\begin{linenomath}\begin{equation}\label{pnga<pnr}
	\pnga < \pnr
\end{equation}\end{linenomath}
holds always (because $p_3>0$). 

(d) The inequalities $\pnga<\pnplus$, $\pnplus<p_2+2p_3$, and $\pnga<p_2+2p_3$ are equivalent and hold if and only if 
\begin{linenomath}\begin{equation}\label{pnplus>pnga}
	p_3>\tfrac19 \;\text{ or }\; p_2 > \sqrt{p_3}-3p_3\,.
\end{equation}\end{linenomath}
Hence, $\pnplus < p_2+2p_3$ is incompatible with $\pnplus \le \pnga$.
If $p_3>\tfrac19$ or if $p_2 > \sqrt{p_3}-p_3$ (which is a restriction only if $p_3<\tfrac19$), then $\pnga < p_2+2p_3$ holds. 

(e) The following relations hold between the bounds presented above:
\begin{linenomath}\begin{equation}
	\sqrt{p_3} -3p_3 < \tfrac12\left(\sqrt{p_3(4+p_3)}-5p_3\right) < \tfrac12\left(\sqrt{p_3(4+p_3)}-3p_3\right) < \sqrt{p_3} -p_3\,.
\end{equation}\end{linenomath}
These are valid for all $0<p_3<1$; equality holds in all cases if $p_3=0$.
\end{remark}

Here is a key lemma:

\begin{lemma}\label{lem:casesF3}
In the region $R$ the following cases can be distinguished:

(1) $p_0>\rFth$ and $\gaFth \mFth <1$ in $R$ if and only if  
\begin{linenomath}\begin{equation}\label{lem_eq1}
	\max\{\pnr,0\} < p_0 < \min\{1-p_2-p_3,p_2+2p_3\}  \,,
\end{equation}\end{linenomath}
where $\pnr < p_2+2p_3$ requires \eqref{pnr<m}.

(2) $p_0=\rFth$ and $\gaFth \mFth < 1$ in $R$ if and only if
\begin{linenomath}\begin{equation}\label{lem_eq2}
	0 < \pnr = p_0 < p_2+2p_3 \,,
\end{equation}\end{linenomath}
where  $0< \pnr < p_2+2p_3$ requires \eqref{pnplus_pos} and \eqref{pnr<m}.

(3) $p_0<\rFth$ and $\gaFth \mFth < 1$ in $R$ if and only if
\begin{linenomath}\begin{equation}\label{lem_eq3}
	\max\{\pnga,0\} < p_0 < \min\{\pnr,p_2+2p_3\} \,,
\end{equation}\end{linenomath}
where $\pnga>0$ if and only if \eqref{pnga>0} holds, and $\pnga<p_2+2p_3$ requires \eqref{pnplus>pnga}.
The following three subcases occur:
\begin{linenomath}\begin{subequations}
\begin{equation}\label{lem_eq3i}
	 \max\{\pnga,0\} < \pnplus < p_0 < \min\{\pnr,p_2+2p_3\}  \,,
\end{equation}\begin{equation}\label{lem_eq3ii}
	\max\{\pnga,0\} < \pnplus = p_0 < \min\{\pnr,p_2+2p_3\}  \,,
\end{equation}
\begin{equation}\label{lem_eq3iii}
	 \max\{\pnga,0\} < p_0 < \pnplus  < \min\{\pnr,p_2+2p_3\}  \,.
\end{equation}
\end{subequations}\end{linenomath}

(4) $p_0<\rFth$ and $\gaFth \mFth = 1$ in $R$ if and only if
\begin{linenomath}\begin{equation}\label{lem_eq4}
	0 < p_0 = \pnga < p_2+2p_3 \,,
\end{equation}\end{linenomath}
where $0 < \pnga < p_2+2p_3$ holds if and only if \eqref{pnga>0} and \eqref{pnplus>pnga} are satisfied.

(5) $p_0<\rFth$ and $\gaFth \mFth > 1$ in $R$ if and only if
\begin{linenomath}\begin{equation}\label{lem_eq5}
	0 < p_0 < \min\{\pnga,p_2+2p_3\} \,.
\end{equation}\end{linenomath}

(6) $p_0\ge\rFth$ and $\gaFth \mFth \ge1$ cannot occur in $R$.

In cases (2) -- (5), $p_0<1-p_2+p_3$ is satisfied if the respective display equation is fulfilled. 
\end{lemma}

The elementary but tedious proof is given in Appendix \ref{sec:Proof_F3}. 
The following theorem characterizes the sign structure of the function $\fFth(x)$ defined in \eqref{def_f}. Figure \ref{fig_F3} illustrates all cases.

\FloatBarrier

\begin{theorem}\label{thm:main_result_F3}
We assume that $(p_0,p_2,p_3)\in R$. 

(1) $\fFth(x)\ge0$ on $[0,1]$ occurs in cases (1) and (2) of Lemma \ref{lem:casesF3}; see Fig.~\ref{fig_F3}A,B.

(2) $\fFth(x)$ changes sign once on $[0,1]$ in case (3) of Lemma \ref{lem:casesF3}. The following subcases occur:

	\quad (i) The sign change occurs below $\PinfFth$ if $p_0>\pnplus$; see Fig.~\ref{fig_F3}C.

	\quad (ii) The sign change occurs at $\PinfFth$ if $p_0=\pnplus$; see Fig.~\ref{fig_F3}D.

	\quad (iii) The sign change occurs above $\PinfFth$ if $p_0<\pnplus$; see Fig.~\ref{fig_F3}E.

(3) $\fFth(x)\le0$ on $[0,1]$ occurs in cases (4) and (5) of Lemma \ref{lem:casesF3}; see Fig.~\ref{fig_F3}F,G. 
\end{theorem}

The proof of this theorem is given in Appendix \ref{sec:Proof_F3}. 
In combination with Proposition \ref{prop:bound_Psurv_gamma} and Lemma \ref{lem:casesF3}, Theorem \ref{thm:main_result_F3} immediately yields the aspired characterization concerning lower and upper bounds for the extinction probabilities $\PnFth$.

\begin{corollary}\label{cor:main_result_F3}
Assume the probability distribution defined in \eqref{pk_F3} with the additional constraint $p_0<p_2+2p_3$. Then the extinction probability by generation $n$, $\PnFth$, has the following properties:

(1) $\PnFth$ satisfies \eqref{PnFLgam} for every $n\ge0$ if and only if 
\begin{linenomath}\begin{equation}\label{lower_bound}
	\pnr \le p_0 < 1-p_2-p_3 \;\text{ and }\; 0<p_0<p_2+2p_3 \,.
\end{equation}\end{linenomath}

(2) $\PnFth$ satisfies \eqref{PnFLgam} for large $n$, and \eqref{PnFLgam_reversed} for small $n$, if and only if
 \begin{linenomath}\begin{equation}\label{no_bound}
	\pnga < p_0 < \pnr  \;\text{ and }\;  0<p_0<p_2+2p_3 \,.
\end{equation}\end{linenomath}

(3) $\PnFth$ satisfies \eqref{PnFLgam_reversed} for every $n\ge0$ if and only if 
\begin{linenomath}\begin{equation}\label{upper_bound}
	 p_0 \le \pnga  \;\text{ and }\;  0<p_0<p_2+2p_3 \,.
\end{equation}\end{linenomath}
\end{corollary}

Analogous statements hold for $\SnFth$, the survival probability until generation $n$ (cf.\ Proposition \ref{prop:bound_Psurv_gamma}). To relate the cases in Corollary \ref{cor:main_result_F3} to each other, it is useful to recall from Remark \ref{rem_props p*}(b) that $\pnplus<\pnr$ if and only if $0<\pnplus$. 

Remark \ref{rem_props p*} informs us that \eqref{pnr<m} is a necessary condition for \eqref{lower_bound} to hold, and $p_2 < \sqrt{p_3} - p_3$ is necessary for \eqref{no_bound} and \eqref{upper_bound}. Furthermore, \eqref{lower_bound} implies $p_2\le\tfrac12$ and $p_0\le\tfrac23$, where $p_0=\tfrac23$ is attained if $p_2=0$ and $p_3=\tfrac13$. Next, \eqref{no_bound} implies $p_2<\tfrac14$ and $p_0<\tfrac13$, where the supremum $\tfrac13$ of $p_0$ is attained at $p_2=0$ and $p_3=\tfrac16$. Finally, \eqref{upper_bound} implies $p_2<\tfrac14$ and $p_0<\tfrac29$, where the supremum $\tfrac29$ of $p_0$ is attained at $p_2=0$ and $p_3=\tfrac19$.

Figure \ref{fig_region_F3} displays the regions defined in statements (1), (2), and (3) of Corollary \ref{cor:main_result_F3}, which are the same as those in (1), (2), and (3) of Theorem \ref{thm:main_result_F3}. The volume of the region defined in (1) is approximately 86.6\% of the total volume of $R$; that of the region in (2) is approximately 10.2\%, and the volume of the region in (3) is approximately 3.2\% of the total volume of $R$.

\begin{remark}\label{rem:p3=0}
The case $p_3=0$ and $p_2>p_0>0$ is treated readily. Retaining the notation from above, we obtain $\mFth=1-p_0+p_2>1$, $\PinfFth=\frac{p_0}{p_2}$, and $\gaFth = 1+p_0-p_2 < 1/\mFth$. Defining $f(x)$ in analogy to \eqref{def_f}, where now $\rFth=\frac{p_0}{1+p_0}$ and $\pFth=\frac{p_2}{1+p_0}$, we obtain $f(x)=\frac{(1-x)(p_0-p_2x)^2}{1+p_0-p_2x}$. Obviously, we have $f(\PinfFth)=f(1)=0$, $f'(\PinfFth)=0$, $f''(\PinfFth)>0$, and it follows immediately that $f(x)\ge0$ on $[0,1]$. Graphs look similar to that in Fig.~\ref{fig_F3}A. In particular, $\PnFth$ satisfies \eqref{PnFLgam} for every $n\ge0$.
\end{remark}

\FloatBarrier

\section{Bounds and approximations for $\Sinfphi$}\label{sec:Bounds for Sinf}
Analytical expressions for the extinction probability $\Pinfphi$, hence for the survival probability $\Sinfphi$, are rarely available. Therefore, the bounds for the extinction and the survival probabilities up to generation $n$ derived on the basis of Proposition \ref{prop:bound_Psurv_gamma} yield little detailed insight. Of course, numerical evaluation is simple and straightforward. Also the bound for the minimum time $T_\varphi(\ep)$ in \eqref{Tep_gen_ineq}, after which survival is `almost' certain in the sense of \eqref{Tep_gen_def}, depends on $\Pinfphi$ and $\gaphi$.
Several authors have derived bounds and approximations for the extinction probability. A classical result is \citetalias{Haldane1927} approximation, who argued that for a Poisson offspring distribution, the probability of survival of a single mutant with a small selective advantage of $s$ is approximately $2s$. 

We start by presenting results of \citet{Quine1976} and \citet{DaleyNarayan1980}, which yield very accurate bounds and approximations for rather general offspring distributions, especially in the slightly supercritical case.  
In Sect.~\ref{subsec:bounds_series}, we assume $\varphi'(1-) = m=1+s$ and that $\varphi$ can be parameterized by $s$ (and other parameters). We derive series expansions of $\Sinfphi$ and $\gaphi$ in terms of $s$ by assuming that $s$ is sufficiently small. Then we highlight the relation to old and recent results on Haldane's approximation for $\Sinfphi$. In Sect.~\ref{subsec:fix_prob_WF} we briefly discuss the relation of survival probabilities in the Galton-Watson process with the diffusion approximation for the fixation probabilities in a finite Wright-Fisher population. In Sects.~\ref{subsec:Poisson_bounds}, \ref{subsec:Bin_bounds}, \ref{subsec:NB_bounds}, and \ref{subsec:GP}, we specify the upper and lower bounds of Quine and of Daley and Narayan, as well as the series expansions of $\Sinfphi$ and $\gaphi$ for the Poisson, the binomial, the negative binomial, and the generalized Poisson distribution, respectively. The accuracy of these bounds and approximations for $\Sinfphi$ is investigated numerically (Table \ref{table_Sinf}). In Sect.~\ref{subsec:GP} we apply the series expansion method to the generalized Poisson distribution and obtain analytical results on the validity of the upper or lower bound, \eqref{SnFLgam} or \eqref{SnFLgam_reversed}, for the time-dependent survival probabilities.

\subsection{Quine's bounds}\label{sec:Quine_bound}
Throughout this and the subsequent sections let  $m=\varphi'(1-)$ denote the mean of $\varphi$, $b=\varphi''(1-)$, $c=\varphi'''(1-)$, and $\si^2 = b+m-m^2$ the variance. In addition to our general assumptions $m>1$ and $0<\si^2<\infty$, we assume $0<c<\infty$. 
We define the quantity
\begin{linenomath}\begin{equation}
	\be := \frac{2(m-1)}{b} \,. 
\end{equation}\end{linenomath}
Under the assumption
\begin{equation}\label{phi_bound_Quine}
	2\be < \min\left\{ 1, \frac{3b}{2c} \right\}\,,
\end{equation}
\citet[][Theorem 2]{Quine1976} derived the following lower and upper bounds for the survival probability $\Sinfphi$:
\begin{equation}\label{Quine_two-sided bounds}
	\LQ_{\varphi} < \Sinfphi < \UQ_{\varphi} \,,
\end{equation}
where
\begin{align}
	\LQ_{\varphi} &= \be + \be^2\, \frac{\varphi'''(1-2\be)}{3b} \,, \\
	\UQ_{\varphi} &= \be + \be^2\, \frac{c}{3b}\Bigl(1-\frac{4c}{3b}\be\Bigr)^{-3/2}\,.
\end{align}
(Quine formulated his result for $\Pinf=1-\Sinf$, and he used $\phi$ instead of $\be$.) 

\subsection{Daley and Narayan's upper bound for $\Sinfphi$}\label{sec:DN_bound}
\citet[][Lemma 3]{DaleyNarayan1980} proved that if
\begin{equation}\label{DN_cond}
	8c(m-1) < 3b^2\,,
\end{equation}
which is equivalent to $2\be < \frac{3b}{2c}$, then 
\begin{equation}
	\Sinfphi < \UDN_{\varphi} \,,
\end{equation}
where
\begin{equation}\label{BDN}
	\UDN_{\varphi} = \frac{3b-3\sqrt{b^2-\tfrac83c(m-1)}}{2c}\,.
\end{equation}
They also showed that condition \eqref{DN_cond} cannot be satisfied if $m\ge 3.2$. In addition, they derived a lower bound.
It is easy to show that $\UDN_{\varphi} < \UQ_{\varphi}$ whenever \eqref{DN_cond} holds (and $m>1$).

For a wide variety of families of probability distributions, \citet{From2007} compared a large number of different upper and lower bounds for $\Sinfphi$ derived by various authors. For $m>1$ and close to 1, he concluded that $\UDN$ is the best upper bound for $\Sinfphi$ among the bounds investigated, and Quine's lower bound is the best lower bound (slightly better then a bound given by \citealt{Narayan1981}). In addition, From derived new, simple, general upper and lower bounds in terms of $p_0$, $p_1$, and $p_2$, which are most useful for large $m$, such as $m>1.5$.

\subsection{Series expansions of $\Sinfphi$ and $\gaphi$}\label{subsec:bounds_series}
In the slightly supercritical case, there is a relatively simple, easily automatized procedure to derive series expansions of these quantities.  
To put this on a firm mathematical basis, we consider families of offspring generating functions $\varphi(x;s)$ depending smoothly on $s\ge0$, typically through various parameters that depend on $s$. We denote partial derivatives of $\varphi(x;s)$ of order ${(k,l)}$ and evaluated at $(x_0,s_0)$ by $\varphi^{(k,l)}(x_0;s_0)$, and we denote $\mu_{kl}=\varphi^{(k,l)}(1;0)$ (we assume that at least the one-sided limits and derivatives at $(1-,0+)$ exist and are finite to the order considered). In particular, we assume that $m(s)= \varphi^{(1,0)}(1-;s)= 1+s$. We define $b=b(s)= \varphi^{(2,0)}(1-;s)$ and $c=c(s)= \varphi^{(3,0)}(1-;s)$. We assume that $\mu_{20}=\lim_{s\to0+}b(s)>0$ and $\mu_{30}=\lim_{s\to0+}c(s)\ge0$. 

To derive a series expansion of $\Sinfphi(s)$ under the assumption that $s$ is small, we set $\Sinfphi(s) = \sum_{i=1}^k  \de_i s^i +O(s^{k+1})$. Then we expand $\varphi(1-\sum_{i=1}^k \de_i s^i;s) - (1-\sum_{i=1}^k \de_i s^i)$ up to $s^{k+1}$ (see Sect.~6.2 in the supplementary \emph{Mathematica} notebook). The resulting coefficient of $s$ vanishes. Equating the coefficients of $s^2, \ldots, s^{k+1}$ to 0, yields $\de_1,\ldots,\de_k$.

With $m=1+s$, the variance is $\si^2 = b-s-s^2$ and we define $\si^2_0=\lim_{s\to0+}\si^2(s)$. Then $\si^2_0 =\lim_{s\to0+}b(s) = \mu_{20}>0$ and we introduce
\begin{equation}\label{theta}
	\theta := \frac{2}{\si^2_0}= \frac{2}{\mu_{20}} \,.
\end{equation}
The approach outlined above yields
\begin{equation}\label{Sinf_series_gen}
	\Sinfphi = \theta s - \de_2 s^2 + \de_3 s^3 + O(s^4) \,,
\end{equation}
where
\begin{subequations}\label{d_2,d_3}
\begin{align}
	\de_2 &= \frac{6\mu_{20}\mu_{21} - 4\mu_{30}}{3\mu_{20}^3} \,, \\  
	\de_3 &= \frac{1}{9\mu_{20}^5}\bigl(18\mu_{20}^2\mu_{21}^2 - 9\mu_{20}^3 \mu_{22} + 16\mu_{30}^2 - 36\mu_{20}\mu_{21}\mu_{30} + 12\mu_{20}^2\mu_{31} - 6\mu_{20} \mu_{40} \bigr)\,.
\end{align}
\end{subequations}
Higher-order terms are readily derived by this method but are increasingly complicated because the coefficient of $s^j$ depends on the mixed partial derivatives of $\varphi(x;s)$ up to order $j+1$ and evaluated at $(1;0)$ (see Sect.~6.2 in the supplementary \emph{Mathematica} notebook).

Given \eqref{Sinf_series_gen}, straightforward calculations yield the following expansion of $\gaphi(s)=\varphi^{(1,0)}(\Pinfphi(s);s)$:
\begin{equation}\label{gaphi_series_gen}
	\gaphi(s) = 1 - s + \ga_2 s^2 -\ga_3 s^3 +  O(s^4)\,,
\end{equation}
where 
\begin{subequations}\label{ga_2,ga_3}
\begin{align}
	\ga_2 &= \frac{2\mu_{30}}{3\mu_{20}^2} \,, \\
	\ga_3 &= \frac{2}{9\mu_{20}^4} \bigl(6\mu_{20}\mu_{21}\mu_{30} - 4\mu_{30}^2 - 3\mu_{20}^2\mu_{31} + 4\mu_{20}\mu_{40} \bigr)\,.
\end{align}
\end{subequations}
Remarkably, the universal coefficient $-1$ of $s$ arises. 
This approximation is useful and provides insight because explicit analytical expressions for $\gaphi$ rarely exist (for a few  exceptions, see below).

We recall from \eqref{mphi*gaphi<1} that $\mphi\gaphi<1$ is equivalent to $\varphi^{(1,0)}(1,s)> \varphiFL^{(1,0)}(1,s)$, which is a necessary condition for $\varphi(x)>\varphiFL(x)$ to hold for $x\in[0,1]$. From \eqref{gaphi_series_gen} we conclude that $\mphi\gaphi<1$ holds if $\ga_2 < 1$ and $s$ is sufficiently small.

The bounds of Quine (Sect.~\ref{sec:Quine_bound}) and of Daley and Narayan (Sect.~\ref{sec:DN_bound}) can be applied to families $\varphi(x;s)$ of pgfs. Interestingly, series expansions of the lower and upper bounds $\LQ_{\varphi}$ and $\UQ_{\varphi}$ in \eqref{Quine_two-sided bounds} and of the upper bound $\UDN_{\varphi}$ in \eqref{BDN} all yield the correct second-order term $\de_2$ in \eqref{d_2,d_3}. The coefficients of $s^3$ differ, and that of $\UDN_{\varphi}$ is closer to the true value $\de_3$ than that of $\UQ_{\varphi}$ (see Sect.~6.3 in the \emph{Mathematica} notebook). From the leading-order term $\be$ of $\LQ_{\varphi}$ and from the expansion \eqref{Sinf_series_gen} (or that of $\UQ_{\varphi}$ or $\UDN_{\varphi}$) we obtain for sufficiently small $s$ the simple bounds
\begin{equation}\label{2s/b<Sinf<th s}
	\be(s) := \frac{2s}{b(s)} < \Sinfphi  < \theta s \,,
\end{equation}
\emph{provided} $\de_2>0$. We have $\de_2>0$ for the Poisson, binomial, and negative binomial distributions, whereas for the generalized Poisson distributed treated below, this holds only for sufficiently small $\la$. A simple example where $\de_2<0$ is the following: let $p_k=0$ for $k\ge3$ and $p_0=\tfrac12-2s$, $p_1=3s$, $p_2=\tfrac12-s$. Then $m=1+s$, $b(s)=1-2s$, and $\Sinfphi=\frac{2s}{1-2s}= 2s + 4s^2 + O(s^3) > 2s = \theta s$.

These bounds and expansions are closely related to the generalized version of \citetalias{Haldane1927} approximation, which in our notation reads
\begin{equation}\label{Haldane_approx_gen}
	\Sinfphi(s) = \theta s + O(s^2) \; \text{ as } \; s\to 0+\,.
\end{equation}
This was derived in a branching-process context, and in various degrees of generality, by  
\citet{Ewens1969}, \citet{Eshel1981}, \citet{Hoppe1992}, and \citet{Athreya1992}; see also \citet[][p.\ 126]{Haccou2005}.  
By contrast, \citet{LessardLadret2007} and \citet{Boenkost_etal2021_weak,Boenkost_etal2021_strong} proved \eqref{Haldane_approx_gen} for certain Markov chain models of Cannings type (a generalization of the Wright-Fisher model), where then the left-hand side is the fixation probability. Indeed, \citet{LessardLadret2007} proved (a generalized version of) $\Sinfphi(s) = \frac{1}{N} + \theta s + o(s)$, where $N$ is fixed as $s\to0$, whence selection is weak relative to random genetic drift. Quite differently, \citet{Boenkost_etal2021_weak,Boenkost_etal2021_strong} assumed that $s$ is asymptotically equivalent to $N^{-b}$ (where $0<b<\tfrac12$ or $\tfrac12<b<1$) as $N\to\infty$. If $0<b<1$, selection is stronger than in the diffusion approximation, where $s$ is asymptotically equivalent to $N^{-1}$.

We note that \eqref{2s/b<Sinf<th s} as well as Quine's bounds in \eqref{Quine_two-sided bounds} imply the generalized version \eqref{Haldane_approx_gen} of Haldane's approximation for Galton-Watson processes because $\lim_{s\to0}b(s)=\mu_{20} = 2/\theta$. 
Interestingly, the lines of research on Haldane's approximation (cited above) and on bounds for the extinction probability and extinction times (e.g., \citealt{Seneta1967}, \citealt{Agresti1974}, \citealt{Quine1976}, \citealt{DaleyNarayan1980}, \citealt{Narayan1981}, \citealt{From2007}) apparently developed independently as no cross references occur; in the latter case, not even to Haldane.

In Table \ref{table_Sinf} we present numerical examples that demonstrate the accuracy of the bounds and approximations presented above. We chose the value $s=0.2$,  despite being high for an advantageous mutant, because for the distributions shown the relative errors vanish rapidly as $s$ decreases below $0.1$.

\subsection{Relation to fixation probabilities in the Wright-Fisher model}\label{subsec:fix_prob_WF}
Following \citet[][p.~120]{Ewens2004}, we define the variance-effective population size by $N_e=N/(\si^2/m)$, where $m$ and $\si^2$ are mean and variance of the offspring distribution. Then the diffusion approximation for the fixation probability of a \emph{single} mutant with selective advantage $s$ in the (haploid) Wright-Fisher model is
\begin{equation}\label{PfixD}
	P^{\rm fix}_D(s,N,N_e) = \frac{1 - e^{-2sN_e/N}}{1 - e^{-2sN_e}}\,.
\end{equation}
If we set $N=1000$, $s=0.1$, $m=1+s$, and $\si^2=\tfrac12(1+s)$, $1+s$, and $5(1+s)$ (so that $N_e=2N, N, \tfrac15N$), then 
$P^{\rm fix}_D \approx 0.3297$, $0.1813$, and $0.0392$, respectively. Interestingly, the survival probabilities in the Galton-Watson process with corresponding fractional linear offspring distributions are $\SinfFL \approx 0.3333$, $0.1818$, and $0.0392$, thus almost identical. With $s=0.1$ we obtain for the Poisson distribution $\SinfPoi \approx 0.1761$, and for the binomial distribution $\SinfBin\approx 0.1763$ (where we set $n=N$). The latter two values are nearly identical to the exact fixation probability $P^{\rm fix}\approx0.1761$ in the standard Wright-Fisher model with $N_e=N$ (computed from the linear system defining the fixed point of the transition matrix; e.g.~\citealt[][p.~87]{Ewens2004}). If $N=N_e=100$ and $s=0.1$, then $P^{\rm fix}\approx0.1758$, $P^{\rm fix}_D\approx0.1813$, $\SinfBin\approx 0.1778$, and $\SinfPoi$ remains unchanged. 

Now we assume $N_e=N$ and the standard Wright-Fisher model. \citet{BuergerEwens1995} proved that the diffusion approximation is always an upper bound for the exact fixation probability $P^{\rm fix}$, and its error is of order $s^2$. In addition, they derived a bound for the relative error of approximations of the form 
\begin{equation}\label{PfixA}
	P^{\rm fix}_A(s,N) = \frac{1 - e^{-A(s)}}{1 - e^{-A(s) N}}\,,
\end{equation}
where $A(s)=a_1 s + a_2s^2$ (in fact, they admitted convergent series). The relative error is of order $s^2$ if $a_1=2$ and $a_2=0$ (yielding the diffusion approximation). They also showed that $a_2$ can be chosen such that the relative error is of order $s^3$. In the haploid case, their equation (4.11) applies and yields $a_2 = -\tfrac43 - \frac{1}{3\nu} + O(e^{-2\nu})$, where $\nu=Ns$ is large, but constant (indeed the coefficient of $e^{-2\nu}$ can be computed explicitly, but is irrelevant in our context). However, this improved, diffusion-like approximation is no longer a global bound for the true $P^{\rm fix}$. Its series expansion in $s$ (with $\nu$ constant) is $2s - \bigl(\frac{8}{3} +\frac{1}{3\nu}\bigr)s^2 + O(s e^{-2\nu})+ O(s^3)$, thus nearly identical to the approximation $\SinfPoi\approx 2s - \frac{8}{3}s^2$ in \eqref{SinfPoi_series} below if $\nu$ is sufficiently large. The diffusion approximation $P^{\rm fix}_D$ has the expansion  $2s - \frac{2}{3}s^2 + O(s e^{-2\nu})+ O(s^3)$. If $s=0.1$, then $P^{\rm fix}_A \approx 0.1758$ if $N=1000$, and $P^{\rm fix}_A \approx 0.1755$ if $N=100$, which are nearly identical to the true values of $0.1761$ and $0.1758$, respectively, in the Wright-Fisher model. B\"urger and Ewens derived also a simple diffusion-like lower bound; it is obtained by setting $A(s)=s/(1+s)$. Its series expansion is $2s-4s^2 + O(s^3)$. It is informative to compare the series expansions of these bounds with those in the following section.

It would be of interest to explore when the survival probability in a Galton-Watson process yields a better approximation for the fixation probability in the Wright-Fisher model with appropriately chosen $N_e$ than the standard diffusion approximation. The work of \citet{LessardLadret2007} and \citet{Boenkost_etal2021_weak,Boenkost_etal2021_strong} (discussed above) could provide a valuable starting point.

\subsection{Poisson distribution}\label{subsec:Poisson_bounds}
For the Poisson distribution, we obtain the following series expansions directly from \eqref{PinfPoi} and \eqref{gaPoi} by using {\it Mathematica} (Sect.~6.4 in the notebook):
\begin{equation}\label{SinfPoi_series}
	\SinfPoi = 2s - \frac{8s^2}{3} + \frac{28s^3}{9} + O(s^4) \,,
\end{equation}
\begin{equation}\label{gaPoi_series}
	\gaPoi = 1 - s + \frac{2s^2}{3} -\frac{4s^3}{9}+ O(s^4)\,.
\end{equation}
These expansions are based on the Taylor series of the Lambert function, $W(x)=\sum_{k=1}^\infty \frac{(-k)^{k-1}}{k!}x^k$, which converges if $|x|<1/e$. The series for $\SinfPoi$ converges if $0\le s < 1$. 

The upper bound of \citet{DaleyNarayan1980} simplifies to  
\begin{equation}
	\UDN_{\rm Poi} = \frac{3 - \sqrt{24/m-15}}{2m} = 2s- \frac{8s^2}{3} + \frac{34s^3}{9} + O(s^4).
\end{equation}
The lower bound  of \citet{Quine1976} becomes
\begin{equation}
	\LQ_{\rm Poi} = 2s - \frac{8s^2}{3} - \frac{10s^3}{3} + O(s^4)\,.
\end{equation}
and the series expansion of the simple lower bound $\be=\frac{2(m-1)}{b}$ is
\begin{equation}
	\be(s) = \frac{2s}{(1+s)^2} = 2s - 4s^2 + 6s^3 + O(s^4)\,.
\end{equation}
The bounds $\UDN_{\rm Poi}$ and $\LQ_{\rm Poi}$ apply if $s<\frac{3}{5}$; $\be$ applies always but becomes very inaccurate if $s\ge0.5$.

\subsection{Binomial distribution}\label{subsec:Bin_bounds}
For the binomial distribution we use the method outlined in Sect.~\ref{subsec:bounds_series} to derive a series expansion of $\SinfBin$. With $m=1 + s$ and $p=\frac{1+s}{n}$, we obtain
\begin{equation}\label{SinfBin_series}
	\SinfBin = \frac{2n}{n-1} s - \frac{4n(2n-1)}{3(n-1)^2} s^2 + \frac{2n(14n^2-17n+5)}{9(n-1)^3} s^3 + O(s^4)
\end{equation}
By differentiation of the generating function, we obtain $\gaBin = \frac{np\PinfBin}{1-p+p\PinfBin}$, which yields after substitution of \eqref{SinfBin_series}:
\begin{equation}\label{gaBin_series}
	\gaBin = 1 - s + \frac{2(n-2)}{3(n-1)}s^2 -\frac{4(n-2)^2}{9(n-1)^2}s^3 + O(s^4)\,.
\end{equation}

The upper bound of \citet{DaleyNarayan1980} becomes
\begin{subequations}
\begin{align}
	\UDN_{\rm Bin} &= \frac{3n\Bigl(1-\sqrt{1-\frac{8(n-2)s}{3(n-1)(1+s)}}\,\Bigr)}{2(n-2)(1+s)}  \\
	&= \frac{2n}{n-1} s - \frac{4n(2n-1)}{3(n-1)^2} s^2 + \frac{2n(17n^2-32n+23)}{9(n-1)^3} s^3 + O(s^4)\,,
\end{align}
\end{subequations}
which is a valid bound if $n\ge2$ and $m\le\tfrac85$. As already noted, the lower bound $\LQ_{\rm Bin}$ of \citet{Quine1976} has the same coefficients of $s$ and $s^2$.

\subsection{Negative binomial distribution}\label{subsec:NB_bounds}
%
For the negative binomial distribution with $m=1+s$ and $p= \frac{r}{r+1+s}$
we obtain by the method outlined in Sect.~\ref{subsec:bounds_series},
\begin{equation}\label{SinfNB_series}
	\SinfNB = \frac{2r}{r+1} s - \frac{4r(2r+1)}{3(r+1)^2} s^2 + \frac{2r(14r^2+17r+5)}{9(r+1)^3} s^3 + O(s^4)
\end{equation}
and 
\begin{equation}\label{gaNB_series}
	\gaNB = 1 - s + \frac{2(r+2)s^2}{3(r+1)} -\frac{4(r+2)^2s^3}{9(r+1)^2} + O(s^4)\,.
\end{equation}

The bound of \citet{DaleyNarayan1980} becomes
\begin{subequations}\begin{align}
	\UDN_{\rm NB} &= \frac{3r}{2(r+2)(1+s)}\biggl(1-\sqrt{1-\frac{8(r+2)s}{3(r+1)(1+s)}} \, \biggr) \\
		&= \frac{2r}{r+1} s - \frac{4r(2r+1)}{3(r+1)^2} s^2 + \frac{2r(17r^2+32r+23)}{9(r+1)^3} s^3 + O(s^4)\,.
\end{align}\end{subequations}
which is a valid bound if $m \le \frac{8(r+2)}{5r+13}$. The simple lower bound $\be$ has the expansion
\begin{equation}
	\be = \frac{2rs}{(r+1)(1+s)^2} = \frac{2r}{r+1} s - \frac{4r}{r+1} s^2 + \frac{6r}{r+1} s^3 + O(s^4)\,.
\end{equation}

\subsection{Generalized Poisson distribution}\label{subsec:GP}
The following generalization of the Poisson distribution was introduced by \citet{ConsulJain1973}:
\begin{linenomath}\begin{equation}
	p_{\rm{GP}}(k) = \frac{\mu(\mu+k\la)^{k-1}}{k!} e^{-\mu-k\la}\,, \quad k=0,1,2,\ldots\,,
\end{equation}\end{linenomath}
where $\mu>0$ and $0\le\la<1$. If $\la=0$, this reduces to the Poisson distribution with $\mu=m$.
\citet[][Chap.~7.2.6]{JohnsonKempKotz2005univariate} call it the Lagrangian Poisson Distribution and summarize relevant properties and relations to other distributions. For a detailed treatment and review of applications consult Chap.~9 of \citet{ConsulFamoye2006}. 
For a relatively simple proof that $\sum_{k=0}^\infty p_{\rm{GP}}(k)=1$, see \citet{Tuenter2000}.

The mean and variance of this unimodal distribution are $m = \dfrac{\mu}{1-\la}$ and $\si^2 = \dfrac{\mu}{(1-\la)^3}$, respectively.
In addition to the coefficient of variation ($\si/m$) also its skew and kurtosis increase to infinity if the mean is held constant and the parameter $\la$ is increased from 0 to 1 \citep[e.g.][Chap.~7.2.6]{JohnsonKempKotz2005univariate}. 
 If $\la>0$, the generating function is given by
\begin{linenomath}\begin{equation}
	\varphiGP(x;\mu,\la) = \exp\Bigl[-\mu\Bigl(1+ \frac{1}{\la}W(-x \la e^{-\la})\Bigr)\Bigr]\,.
\end{equation}\end{linenomath}

We demonstrate the utility of our series expansion by applying it to this distribution. Otherwise, it is difficult to analyze in our context because, apparently, the survival probability $\SinfGP$ cannot be expressed in terms of known functions. All calculations, algebraic and numeric, can be found in detail in Sect.~7 of the supplementary {\it Mathematica} notebook.

If $m=1+s$, then $\mu=(1+s)(1-\la)$, $b_0=b(0)=\frac{1}{(1-\la)^2}$, $c_0=c(0)=\frac{1+2\la}{(1-\la)^4}$, $\mu_{21}=1+b_0$, $\mu_{31}=\frac{6}{(1-\la)^2}$, $\mu_{40}= \frac{1+\la(6+9\la-\la^3)}{(1-\la)^6}$, and $\theta = 2/b_0 = 2(1-\la)^2$.
Therefore, \eqref{Sinf_series_gen} and \eqref{d_2,d_3} yield
\begin{align}\label{SinfGP_series}
	\SinfGP &= 2(1-\la)^2s - \frac{2}{3}(1-\la)^2 (4-10\la +3\la^2)s^2 \notag \\
	&\quad + \frac{4}{9}(1-\la)^3(7-31\la+21\la^2-3\la^3)s^3 + O(s^4)\,.
\end{align}
The coefficient of $s^2$ is positive if and only if $\la < \frac{5-\sqrt{13}}{3}\approx 0.4648$. Therefore, $\SinfGP > \theta s$ for small $s$ if $\la > \frac{5-\sqrt{13}}{3}$. The coefficient of $s^3$ is positive if $\la\lessapprox0.2750$.
From \eqref{gaphi_series_gen} and \eqref{ga_2,ga_3} we obtain
\begin{equation}\label{gammaGP_series}
	\gaGP = 1-s + \frac{2}{3}(1+2\la)s^2 - \frac{4}{9}(1+7\la+\la^2)s^3+ O(s^4)\,.
\end{equation}
We do not present the bounds of Quine and of Daley and Narayan because the expressions are quite complicated. However, in Table~\ref{table_Sinf} numerical values are shown.

\begin{table}[t]
\centering
\begin{tabular}{l|c|c|cccc|c}
\hline\hline
	\multirow{2}*{$s=0.2$}& $\varphiBin$ & $\varphiNB$ & \multicolumn{4}{|c|}{$\varphiGP$:\quad $\la=$}& $\varphiFL$ \\	
					& $n=5$ & $r=5$ & 0 &  0.2 &  0.5 & 0.9 & $\pi=0.2$ \\
	\hline
		  $\be$  & 0.3472 & 0.2315 & 0.2778 & 0.1891 & 0.0794 & 0.00333 & 0.6667 \\
		  $\LQ_\varphi$   & 0.3673 & 0.2444 & 0.2936 & 0.1993 & 0.0832 & 0.00346 & 0.7018\\
		  $\Sinfphi$ & 0.3804 & 0.2668 & 0.3137 & 0.2228 & 0.1003 & 0.00466 & 0.8000 \\
$\Sinf_{\rm{ser}}$& 0.3875 & 0.2670 & 0.3182 & 0.2228 & 0.1001 & 0.00466 & 0.8000 \\
		  $\UDN_\varphi$  & 0.3823 & 0.2733 & 0.3183 & 0.2317 & 0.1158 & --     & 0.8453\\
	    $\theta s$& 0.5000 & 0.3333 & 0.4000 & 0.2560 & 0.1000 & 0.00400 & 0.8000 \\
	\hline\hline
\end{tabular}
\caption{The table shows values of $\Sinfphi$ and its bounds and approximations for $s=0.2$. The data confirm the analytical results of \citet{Quine1976} and \citet{DaleyNarayan1980} that $\be \le \LQ_\varphi \le \Sinfphi \le \UDN_\varphi$. $\Sinf_{\rm{ser}}$ denotes the series expansion up to order $s^3$ given in \eqref{Sinf_series_gen}. Its relative error to $\Sinfphi$ is always smaller than that of $\LQ_\varphi$ and also than that of $\UDN_\varphi$ except for $\varphiBin$. Note that $\Sinf_{\rm{ser}}<\SinfGP$ if $\la=0.5$ and $\la=0.9$ because then the coefficient $\de_2$ is negative; see the text below \eqref{SinfGP_series}. For $\varphiGP$ with $\la=0.9$, $\UDN_{\rm GP}$ yields a complex value because condition \eqref{DN_cond} is violated. The last column contains an example in which the variance of the offspring distribution is very small, so that $\theta$ is large ($\theta=4$). For simplicity, we chose a fractional linear distribution, for which $\SinfFL=\frac{1-\pi}{\pi}s$. In order to achieve $m=1+s$, we chose $\rh=\pi(1+s)-s$; see Sect.~\ref{sec:FL}. Also in this case, the approximations are quite accurate, even if not needed for this distribution. The last line shows the generalized version of Haldane's approximation.
}\label{table_Sinf}
\end{table}

\begin{figure}[t!]
\vspace{0mm}
\centering
\includegraphics[width=0.9\textwidth]{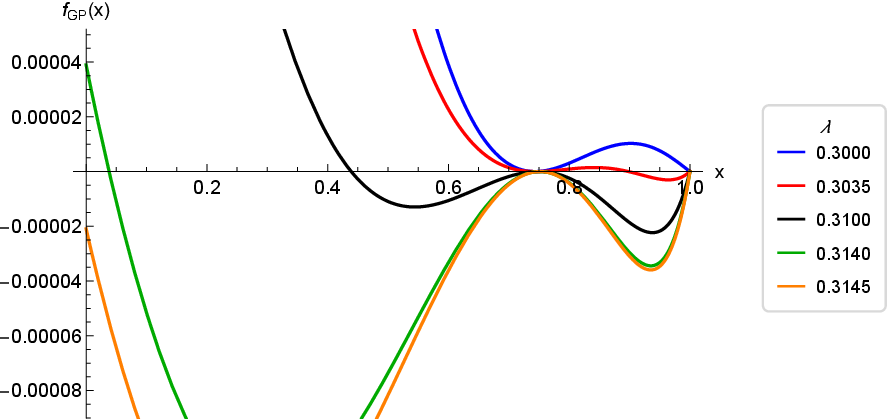} 
\caption{Possible shapes of graphs of $\fGP(x;s,\la)$. 
We chose $s=0.3$ for good visibility. Then $\PinfGP\approx0.7435$ if $\la=0.30$, and $\PinfGP\approx0.7515$ if $\la=0.3145$. At the critical value $\la_{c_1}\approx0.30160$, $\fGP'(1)$ changes sign;  at $\la_{c_2}\approx0.30596$, $\fGP''(\PinfGP)$ changes sign; at $\la_{c_0}\approx0.31433$, $\fGP(0)$ changes sign.
} 
\label{fig_F2}
\end{figure}

Using the series expansions of $\SinfGP$ and $\gaGP$, approximations for the parameters $\pi_{\rm{GP}}$ and $\rh_{\rm{GP}}$ of the bounding fractional linear pgf, $\varphiGP(x)$, can be computed. As in previous sections, we define
\begin{equation}
	\fGP(x_;s,\la) = \varphiGP(x;(1+s)(1-\la),\la) - \varphiFL(x;\pi_{\rm{GP}},\rh_{\rm{GP}})\,,
\end{equation}
where usually we omit the dependence on $s$ and $\la$.
The exact version requires numerical evaluation of $\SinfGP$, $\gaGP$, $\pi_{\rm{GP}}$, and $\rh_{\rm{GP}}$. The qualitatively different shapes are shown in Fig.~\ref{fig_F2}.

The analytical results below are based on calculating $\fGP(x)$ by employing the  series expansions in \eqref{SinfGP_series} and \eqref{gammaGP_series}. We obtain 
\begin{align}\label{fGP[0]}
	\fGP(0) &= e^{-1+\la} -\frac{1}{3-4\la+2\la^2} \notag \\
		&\quad- \Bigl(e^{-1+\la}(1-\la) - \frac{2(1-\la)^2(4+2\la+3\la^2)}{3(3-4\la+2\la^2)^2}\Bigr)s + O(s^2)\,.
\end{align}
By series expansion of $\la$ around the value at which the term of order $1$ vanishes, we find that $\fGP(0)>0$ if and only if $\la < \la_{c_0}$, where
\begin{equation}
	\la_{c_0} \approx0.25915 + 0.1997s 
\end{equation}
provides an accurate approximation if $s\lessapprox0.3$. For instance, if $s=0.1$, the approximation yields $\approx0.27912$ and the numerically determined exact value is $\la_{c_0}\approx0.27857$. 

We recall that by our construction of $\varphiFL(x;\pi_{\rm{GP}},\rh_{\rm{GP}})$, we have $\fGP'(1)= 1+s - \gaGP^{-1}$; cf.~\eqref{mphi*gaphi<1}. Because $(1+s)\gaGP = 1 - \tfrac13(1-4\la)s^2 +\tfrac29(1-8\la-2\la^2)s^3 + O(s^4)$, we find that $\fGP'(1)<0$ if and only if  
$\la<\la_{c_1}$, where 
\begin{equation}
	\la_{c_1}\approx\tfrac14\bigl(1+\tfrac34s\bigr)\,.
\end{equation}
If $s=0.1$, then the numerically precise value is $\la_{c_1}\approx0.26820$, and the simple approximation yields $0.26875$.

For the second derivative of $\fGP$ at $\PinfGP$ we obtain
\begin{equation}
	\fGP''(\PinfGP) = \frac{1-4\la}{3(1-\la)^2}\,s - \frac{1+4\la-74\la^2+12\la^3+3\la^4}{9(1-\la)^2}\,s^2 + O(s^3) \,.
\end{equation}
The first term is positive if $\la<\tfrac14$ and $1+4\la-74\la^2+12\la^3+3\la^4>0$ if $\la<0.14867$. A rough approximation for the critical value $\la_{c_2}$ of $\la$, below which $\fGP''(\PinfGP)>0$ holds, is 
\begin{equation}
	\la_{c_2} \approx \tfrac14 + 0.202s\,.
\end{equation}
If $s=0.1$, then the numerically precise value is $\la_{c_2}\approx0.26967$.

Clearly, $\fGP(x)\ge0$ can hold for $x\in[0,\PinfGP]$ only if $\fGP(0)\ge0$ and 
$\fGP''(\PinfGP)\ge 0$. Both inequalities are satisfied if $0\le\la \le \la_{c_2}$ because our approximations satisfy $\la_{c_2}<\la_{c_0}$ (if $s<4$). 

If $\la_{c_2}<\la<\la_{c_0}$, then $\fGP''(\PinfGP) < 0$ (hence $\fGP(x)<0$ close to $\PinfGP$), and $\fGP(0)>0$. Therefore, $\fGP(x)$ changes sign between $x=0$ and $x=\PinfGP$. Finally, if $\la>\la_{c_0}$ then $\fGP(x)<0$ near $x=0$ and near $x=\PinfGP$. Numerical results suggest that in this case $\fGP(x)<0$ between $0$ and $\PinfGP$, whereas $\fGP(x)>0$ if $0\le \la <\la_{c_2}$. Thus, by Proposition \ref{prop:bound_Psurv_gamma}, the inequalities \eqref{PnFLgam} and \eqref{SnFLgam} hold if $0\le \la <\la_{c_2}$, and the opposite inequalities hold if $\la>\la_{c_0}$. 
The qualitatively different cases are shown in Fig.~\ref{fig_F2} for $s=0.3$. Then the range of values $\la$ in which $\fGP(x)$ changes sign between 0 and $\PinfGP$ is approximately $(0.30596,0.31433)$. If $s\approx0$, then this interval is approximately $(0.25,0.25915$). 

Based on these and additional numerical results (not shown), we conjecture that the inequalities \eqref{PnFLgam} and \eqref{SnFLgam} hold if $0\le \la <\la_{c_2}$, and the opposite inequalities hold if $\la>\la_{c_0}$. In the region $(\la_{c_2},\la_{c_0})$, $\PnGP$ satisfies \eqref{PnFLgam} and \eqref{SnFLgam} for small $n$, and the opposite inequalities for large $n$. This differs from the finitely supported distributions studied in Sect.~\ref{sec:F3}; cf.~Corollary \ref{cor:main_result_F3}(2) and Fig.~\ref{fig_F3}.

\FloatBarrier

\section{Applications}\label{sec:applications}

We begin by investigating the accuracy of the approximations for $T_\varphi(\ep)$ and $\Snphi$ which play a key role in our major application in the final Section \ref{subsec:spread_favorable_mutant}.

\subsection{Convergence time $T_\varphi(\ep)$ of survival probabilities $\Sn$}\label{subsect:Tep}
First, we apply our results to $T_\varphi(\ep)$, the number of generations until the survival probability $\Sn$ differs from the eventual survival probability $\Sinf$ by a factor of at most $1+\ep$; see \eqref{Tep_gen_def}. 
We define 
\begin{linenomath}\begin{equation}\label{Tep_gen_app}
	T_{\rm{app}}(\ep) = \Biggl\lceil\frac{\ln\bigl(\bigl(1+\frac{1}{\ep} \bigr) \Pinfphi \bigr)}{-\ln \gaphi} \Biggr\rceil\,,
\end{equation}\end{linenomath}
where $\lceil z \rceil$ denotes the least integer greater than or equal to $z$.
According to \eqref{Tep_gen_ineq}, this is an upper bound for the true $T_{\varphi}(\ep)$ if \eqref{varphiFL<varphi} holds. It is a lower bound if the reversed inequality \eqref{varphiFL>varphi} holds, and it will serve as an approximation if none of the two holds.

Throughout this section we set $m=1+s$.
By using $\Pinfphi\approx 1-\theta s + \de_2 s^2$ and $\gaphi\approx 1-s + \ga_2s^2 - \ga_3s^3$, where $\theta$, $\de_2$, $\ga_2$ and $\ga_3$ are given in \eqref{theta}, \eqref{d_2,d_3} and \eqref{ga_2,ga_3}, respectively, we obtain from \eqref{Tep_gen_app} by series expansion in $s$
\begin{linenomath}\begin{align}\label{Tep_gen_series_s1}
	T_{\rm{app}}(\ep) &= \Bigl(\frac{1}{s}- \frac12 + \ga_2 \Bigr)\ln\bigl( 1+\frac{1}{\ep}\bigr) - \theta
	\notag \\
	&\quad + \left( \de_2 + \frac12\theta(1-2\ga_2-\theta) + \bigl(\ga_2^2-\ga_3 - \frac1{12} \bigr) \ln\bigl(1 + \frac{1}{\ep}\bigr)  \right)s + O(s^2) \,.
\end{align}\end{linenomath}
Here, $\ln\bigl( 1+\frac{1}{\ep}\bigr)\approx 2.4, 4.6, 6.9$ if $\ep=0.1, 0.01, 0.001$, respectively.
We note that the terms in \eqref{Tep_gen_series_s1} remain unchanged under expansions of $\Pinfphi$ and $\gaphi$ to arbitrary order. 
Now we define
\begin{linenomath}\begin{equation}\label{Tep_gen_series_s0}
	T_{\rm{ser}}(\ep) = \left\lceil\Bigl(\frac{1}{s}- \frac12 + \ga_2 \Bigr)\ln\bigl( 1+\frac{1}{\ep}\bigr) - \theta \right\rceil\,.
\end{equation}\end{linenomath} 
Then $T_{\rm{ser}}(\ep)\ge T_{\rm{app}}(\ep)$ for sufficiently small $s$ if the coefficient of $s$
 in \eqref{Tep_gen_series_s1} is negative. 

%

\begin{table}[th!]
\centering
\begin{tabular}{ll|c|c|ccccc|c}
\hline\hline
	&& $\varphiBin$ & $\varphiNB$ & \multicolumn{5}{|c|}{$\varphiGP$:\quad $\la=$}&\multirow{2}*{$\tfrac{1}{s}\ln(1+\tfrac{1}{\ep})$} \\	
	$s=0.01$ && $n=5$ & $r=5$ & 0 &  0.1 &  0.259 &  0.5 & 0.9 & \\
	\hline
	\multirow{3}{4em}{$\ep = 0.01$}
	& $T_\varphi(\ep)$  & 458 & 461 & 459 & 461 & 463 & 467 & 474 &\multirow{3}*{462} \\
	& $T_{\rm{app}}(\ep)$ & 460 & 462 & 461 & 462 & 463 & 465 & 468 & \\
	& $T_{\rm{ser}}(\ep)$ & 460 & 462 & 461 & 462 & 463 & 465 & 468 & \\
	\hline\hline
	$s=0.1$ &&  &  & 0 &  0.1 &  0.276 &  0.5 & 0.9 &  \\
	\hline
	\multirow{3}{4em}{$\ep = 0.1$} 
	& $T_\varphi(\ep)$  & 21 & 23 & 22 & 23 & 25 & 27 & 31 &\multirow{3}*{24} \\
	& $T_{\rm{app}}(\ep)$ & 22 & 23 & 23 & 24 & 25 & 26 & 28 & \\
	& $T_{\rm{ser}}(\ep)$ & 22 & 23 & 23 & 24 & 25 & 26 & 28 & \\
	\hline
	\multirow{3}{4em}{$\ep = 0.01$}
	& $T_\varphi(\ep)$  & 43 & 46 & 45 & 46 & 48 & 51 & 56 &\multirow{3}*{47} \\
	& $T_{\rm{app}}(\ep)$ & 44 & 46 & 45 & 46 & 48 & 50 & 53 &  \\
	& $T_{\rm{ser}}(\ep)$ & 44 & 46 & 45 & 46 & 48 & 50 & 53 &  \\ 
	\hline
	\multirow{3}{4em}{$\ep = 10^{-4}$}
	& $T_\varphi(\ep)$  & 89 & 93 & 91 & 93 & 96 & 100 & 108 &\multirow{3}*{93} \\
	& $T_{\rm{app}}(\ep)$ & 90 & 93 & 92 & 94 & 96 & 99 & 105 &\\
	& $T_{\rm{ser}}(\ep)$ & 90 & 93 & 92 & 94 & 96 & 100 & 105 &  \\
	\hline\hline
	$s=0.3$ &&  &  & 0 &  0.2 &  0.312 &  0.5 & 0.9 &  \\
	\hline
	\multirow{3}{4em}{$\ep = 0.01$}
	& $T_\varphi(\ep)$  & 13 & 15 & 14 & 16 & 17 & 19 & 24 &\multirow{3}*{16}\\
	& $T_{\rm{app}}(\ep)$ & 13 & 15 & 15 & 16 & 17 & 19 & 22 & \\
	& $T_{\rm{ser}}(\ep)$ & 13 & 15 & 15 & 16 & 17 & 18 & 20 & \\
	\hline\hline
\end{tabular}
\caption{The table shows values of $T_\varphi(\ep)$ for $\varphi=\varphiBin$, $\varphiNB$, and $\varphiGP$ with $m=1+s$ and $s$, $\ep$, and the other parameters as indicated. Here, $T_\varphi(\ep)$ is the exact time defined in \eqref{Tep_definition_FL} and computed by iterating the generating function $\varphi$. $T_{\rm{app}}(\ep)$ is the approximation defined in \eqref{Tep_gen_app} and computed from the numerically exact values of $\Pinfphi$ and $\gaphi$. The series approximation $T_{\rm{ser}}(\ep)$ is defined in \eqref{Tep_gen_series_s0}. The final column contains the values for the simple approximation shown on its top. }
\label{table_Tep}
\end{table}

We recall from \eqref{SinfBin_series} and \eqref{gaBin_series} that for the binomial distribution $\theta=\frac{2n}{n-1}$ and $\ga_2=\frac{2(n-2)}{3(n-1)}$; from \eqref{SinfNB_series} and \eqref{gaNB_series} that for the negative binomial distribution $\theta=\frac{2n}{n+1}$ and $\ga_2=\frac{2(n+2)}{3(n+1)}$; and from \eqref{SinfGP_series} and \eqref{gammaGP_series} that for the generalized Poisson distribution $\theta=2(1-\la)^2$ and $\ga_2=\frac{2}{3}(1+2\la)$.

In Table~\ref{table_Tep} we compare exact values of $T_{\varphi}(\ep)$, obtained by iteration of the generating function, with the approximation $T_{\rm{app}}(\ep)$ in \eqref{Tep_gen_app} and its simple series approximation $T_{\rm{ser}}(\ep)$. We show results for the binomial and negative binomial distribution, the generalized Poisson distribution, and the simple approximation $\frac{1}{s}\ln(1+\frac{1}{\ep})$. We note that if $n\ge10$ and $r\ge10$, the values $T_{\varphi}(\ep)$ for the binomial and the negative binomial are (nearly) identical to those for the Poisson distribution ($\la=0$). For the generalized Poisson distribution, for each $s$ the middle value of $\la$ is chosen from the small range for which $\fGP(x)$ changes sign between 0 and $\PinfGP$ (see Sect.~\ref{subsec:GP}). For the two smaller values of $\la$, $\fGP(x)\ge0$ for every $x$ (and \eqref{SnFLgam} holds), and for the two larger values $\fGP(x)\le0$ (and the opposite of \eqref{SnFLgam} holds). For the middle value of $\la$, the maximum of $|\fGP(x)|$ is extremely small (e.g.\ Fig.~\ref{fig_SnGP}), thus the fractional linear approximation is extremely accurate. 

The data in the table show that mostly the deviations from the very simple approximation $\frac{1}{s}\ln(1+\frac{1}{\ep})$ are small. Larger deviations than shown here occur, for instance, if the variance of the offspring distribution is very small, because then $\theta$ becomes large and $T_\varphi(\ep)$ is reduced. Also a high skew contributes to larger deviations, as is visible for the generalized Poisson distribution with large $\la$, because then $T_\varphi(\ep)$ is increased.

\subsection{Survival probability up to generation $n$}\label{sec:appl_ext_times}
Here, we study the accuracy of the bounds and approximations for the survival probability $\Snphi$ up to generation $n$ and, equivalently, for the extinction probability $\Pnphi$ by generation $n$.
From Proposition \ref{prop:bound_Psurv_gamma} and eq.~\eqref{PnFLgam} we obtain 
\begin{equation}
	 \frac{1-\gaphi^n}{1-\gaphi^n \Pinfphi} \le \frac{\Pexn_\varphi}{\Pinfphi} \le 1
\end{equation}
and, by expansion of the left-hand side, the approximation
\begin{equation}\label{Pexn/Pinf_gen_series}
	 \frac{\Pexn_\varphi}{\Pinfphi} \approx 1 - \frac{\theta}{n+\theta} + \frac{n\bigl(\theta(n+1)+2\de_2-2\theta\ga_2\bigr)}{2(n+\theta)^2}\,s + O(s^2)\,.
\end{equation}
This is accurate if, approximately, $sn<1$.

For the generalized Poisson distribution, \eqref{Pexn/Pinf_gen_series} becomes 
\begin{equation}
	\frac{\Pexn_{\rm{GP}}}{\PinfGP} \approx 1 -\frac{2(1-\la)^2}{n+2(1-\la)^2}+ \frac{(1-\la)^2\bigl(1+\frac{1}{3n}(7-28\la + 6\la^2)\bigr)}{\bigl(1+\frac{2}{n}(1-\la)^2\bigr)^2}\,s + O(s^2)\,.
\end{equation}
Series expansions of $\Snphi/\Sinfphi$ about $s=0$ are less informative because both terms converge to 0 as $s\to 0$.

\begin{figure}[t!]
\vspace{0mm}
\centering
\includegraphics[width=0.9\textwidth]{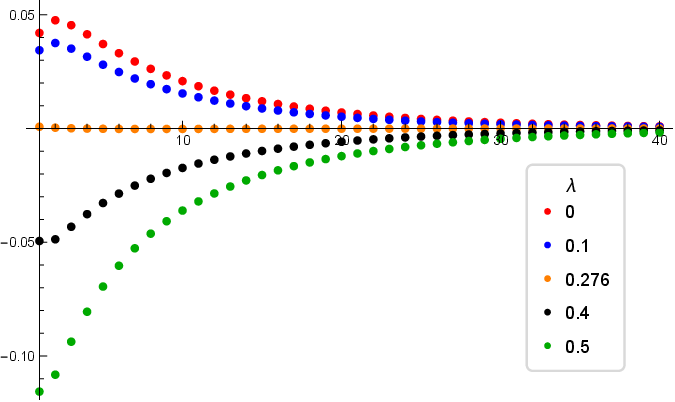} 
\caption{Relative errors of survival probabilities by generation $n$, $(\Sn_{\rm{app}}-\SnGP)/\SnGP$, for the generalized Poisson distribution. In all cases, $s=0.1$. The values of $\la$ are given in the legend; $\la=0$ yields the Poisson distribution. If $\la=0.276$, then $\Sn_{\rm{app}}-\SnGP$ changes sign between $n=3$ and $n=4$; cf.\ Fig.~\ref{fig_F2}. We note that for given $m=1+s$, $\varphiGP$ and $\varphiNB$ have the same variance if $\la=1-\sqrt{r/(r+1+s)}$. With $r=5$ this yields $\la\approx0.095$, $\SinfGP\approx0.14841$, $\SinfNB\approx0.14834$. Thus, on this scale of resolution, the blue curve would be almost indistinguishable from the corresponding curve for the negative binomial with $r=5$.  
} 
\label{fig_SnGP}
\end{figure}

We define
\begin{equation}
	\Sn_{\rm{app}} = \frac{\Sinf_\varphi}{1-\gaphi^n(1-\Sinf_\varphi)}\,,
\end{equation}
which is the right-hand side of \eqref{SnFLgam}. Because of its versatility, we use the generalized Poisson distribution to illustrate the accuracy of $\Sn_{\rm{app}}$ if taken as approximation. From the results in Sect.~\ref{subsec:GP} and Proposition \ref{prop:bound_Psurv_gamma}, we expect that $\Sn_{\rm{app}}$ is an upper bound for $\SnGP$ if $\la<\la_{c_2}$, and a lower bound if $\la>\la_{c_0}$, where $\la_{c_0}>\la_{c_2}$.
Figure~\ref{fig_SnGP}, which displays the relative error $(\Sn_{\rm{app}}-\SnGP)/\SnGP$ for several values of $\la$, confirms this. Table \ref{table_Tep} informs us that for the parameters shown in Fig.~\ref{fig_SnGP}, it takes between $22$ and $31$ generations for $\SnGP$ to decay below $1.1\SinfGP$.

\subsection{The spread of a favorable mutant in a finite population}\label{subsec:spread_favorable_mutant}
Branching process methods have been applied since the early times of population genetics to study the survival of new mutants \citep{Fisher1922,Haldane1927}. Generalizations of Haldane's approximation for the fixation probability of an advantageous mutant are discussed in Sect.~\ref{subsec:bounds_series} along with the relevant references. Another, more recent branch of research has been concerned with  the evolution of the distribution of a favorable mutant in a finite population of constant size $N$. Despite the great utility of diffusion-approximation methods to study the probability of and the expected time to fixation of a mutant \citep[e.g.][]{Charlesworth2020}, they are apparently not conducive to study the time course of the distribution of allele frequencies. 

\citet{DesaiFisher2007} and \citet{UeckerHermisson2011} conditioned on fixation of the favorable mutant and employed branching process methods to approximate its evolution by a deterministic increase starting with the random variable $W=\lim_{n\to\infty}Z_n/m^n$ \citep[e.g.][p.~154]{Haccou2005} as initial condition, whose absolutely continuous part ($W^+$) describes the stochasticity that accumulates during the spread of the mutant. \citet{MartinLambert2015} employed a variant of this approach and approximated the initial and the final phase by a Feller process conditioned on fixation. They derived a semi-deterministic approximation for the distribution of the (beneficial) allele frequency at any given time. 

\citet{GB2024}, in the following abbreviated GB2024, developed a related approach that conditions on survival up to the current generation. This has the advantage that the mean number of mutants is described correctly for the initial generations, and also the variance is approximated very closely. They described the initial phase by a supercritical Galton-Watson process and combined it with the deterministic diallelic selection equation in such a way that the relative frequency of mutants in generation $n$, $X_n$, is given by $X_n=Y_n/(Y_n+N)$, where $Y_n$ is the random variable with the exponential distribution $\Psi_n(y)=1-e^{-\la_n y}$, and $\la_n=\Sn/m^n$ (where the subscript $\varphi$ is suppressed in this section). The underlying rational was to approximate the distribution of the absolutely continuous part $W_n^+$ of $Z_n/m^n$ by an exponential distribution because the limiting distribution $W^+$ is often nearly exponential. It is exponential with parameter $\la=\Sinf$ for  fractional linear offspring distributions, and also $W_n^+$ can be imbedded into the exponential distribution with parameter $\Sn$ (Appendix A in GB2024). Thus, $\Psi_n$ approximates the distribution of $Z_n$ conditioned on $Z_n>0$. 

The discrete distribution of the mutant in generation $n$ (which started as a single copy in generation $0$) can be approximated by the density
\begin{subequations}\label{g(p)}
\begin{linenomath}\begin{equation} 
	g_{a_n} (x) = \frac{a_n}{(1 - x)^2} \cdot  \exp \left[ - a_n \frac{x}{1-x}\right]\,,
\end{equation}\end{linenomath}
where
\begin{linenomath}\begin{equation}\label{a(t)}
	a_n = a_n(m, N, \Sn)  = N \Sn/m^n \,.
\end{equation}\end{linenomath}
\end{subequations}
The structure of $g_{a_n}(x)$ is the same as that of the density $\be_t$ in \citet{MartinLambert2015}, except that their  $\be_t$ (corresponding to our $a_n$) decays exponentially with $t$ and a constant parameter (in our notation $2s$), whereas our $a_n$ has the additional dependence on $\Sn$. 
For large $N$, such as $N\ge1000$, the density in \eqref{g(p)} provides a very accurate approximation for the allele frequency distribution in the corresponding Wright-Fisher model. For the initial generations, it is much more accurate than previous approximations that conditioned on eventual fixation of the mutant (see Sect.~3.3 in GB2024).

One of the main applications in GB2024 of this result was the derivation of explicit formulas for the time dependence of the mean and the genetic variance of a quantitative trait under exponential directional selection (their Propositions 4.3 and 4.11). They assumed that the trait is determined by an underlying infinite-sites model, i.e., every new mutation that contributes to the trait occurs at a new locus ($=$ site), so that many mutants can segregate simultaneously. A key \emph{assumption} was that the offspring distribution is such that the survival probabilities $\Sn$ can be bounded above as in \eqref{SnFLgam}. In Sects.~\ref{sec:Poisson} -- \ref{sec:F3}, we \emph{proved} that Poisson, binomial and some other distributions indeed satisfy this condition. Below we outline for an important special case why this bound for $\Sn$ is essential for the proofs of the results in GB2024 on the evolution of a quantitative trait. In GB2024 this was buried beneath technical complications (see their Appendix D.4).

Given \eqref{g(p)}, the within-population variance of the distribution of mutants in generation $n$ at a single locus is
\begin{linenomath}\begin{equation}\label{within-var}
	\ga(n) = \int_0^1 x(1-x)g_{a_n}(x)\,dx = {a_n}(1+{a_n})e^{a_n}E_1({a_n}) - {a_n} \,,
\end{equation}\end{linenomath}
where $E_1(a)=\int_a^\infty x^{-1}e^{-x}\,dx$ denotes the exponential integral. The model in GB2024 assumes that new mutations occur according to a Poisson process, each mutation at a new site. In the simplest case, which we assume here to avoid plenty of technical detail, each mutation contributes the same effect $\a>0$ to the trait, and its fitness (expected number of offspring) is $m = e^{s\a}$, $s>0$. It was shown that the variance of the trait at time $\tau$ is 
\begin{linenomath}\begin{equation}
	V_G(\tau)  = \Th \a^2 \int_0^\tau S^{([t])} \ga([t])\, dt \,,	\label{VGtau}
\end{equation}\end{linenomath}
where $\Th$ is the expected number of mutations occurring per time unit (generation) in the total population, and $[t]$ denotes the nearest integer to $t$. Because in GB2024 a distribution of mutation effects $\a$ was admitted, an additional integration with respect to $\a$ occurs in eq.~(4.15) of GB2024, which yields \eqref{VGtau} for equal effects. 

In the absence of an explicit formula for $\Sn$ the integrand in \eqref{VGtau} needs to be computed recursively. An explicit formula is available only for fractional linear distributions. Using the present results on bounding the survival probabilities as in \eqref{SnFLgam}, we can quantify when $\Sn$ is sufficiently close to $\Sinf$ so that the integration in \eqref{VGtau} can be simplified by approximating $S^{([t])}$ by $\Sinf$ for sufficiently large $t$. Indeed, in this case the integral can be calculated explicitly and the error terms can be derived. We focus on the important limiting case $\tau\to\infty$, i.e., when the per-generation response of the trait mean and the expected variance become constant due to the balance of loss and fixation of new recurrent mutations. For this case, the basic ideas can be presented without excessive technical detail. 

Following Proposition 4.5 of GB2024, we consider $V_G^\infty=\lim_{\tau\to\infty}V_G(\tau)$ and define $\ga^\infty(t)$ analogously to $\ga(n)$ but with $a^\infty(t)=N\Sinf/m^t=N\Sinf e^{-s\a t}$ instead of $a_n$. We define
\begin{equation}
	V_1^\infty = \int_0^\infty \ga^\infty(t)\, dt = \frac{1}{s\a} N\Sinf e^{N\Sinf} E_1(N\Sinf)\,.
\end{equation}
Then $\a^2V_1^\infty$ is the total variance contributed to the trait by a single mutant during its sweep to fixation (conditioned on its fixation).

GB2024 imposed the assumptions that (i)
\begin{equation}\label{assump1}
	N s^K = C^K \; \text{ as } \; N\to\infty\,,
\end{equation}
where $C>0$ is an arbitrary constant and the constant $K$ satisfies $K>2$, and (ii) the offspring distribution satisfies \eqref{SnFLgam}.
For equal effects $\a$, GB2024 obtained from their Proposition 4.5,
\begin{equation}\label{VGinfty}
	V_G^\infty = \Th \Sinf \a^2 V_1^\infty + O(N^{-K_1/K})\,,
\end{equation}
where $K>K_1+1$ and $K_1>1$ is an arbitrary constant (their Remark 4.6). 
Equation \eqref{VGinfty} shows that to leading order in $N$, the asymptotic variance of the trait depends only on the contribution of mutations that become fixed. The contribution of mutations that are lost is of smaller order. 

By recalling \eqref{Sinf_series_gen} and that $s$ there corresponds to $e^{s\a}-1\approx s\a$ here, we obtain $\Sinf\approx \theta s\a - \de_2(s\a)^2$. A well known asymptotic expansion of $e^xE_1(x)$  yields $V_1^\infty \approx \frac{1}{s\a}\bigl(1-\frac{1}{N\Sinf}\bigr)$ if $N\Sinf$ is sufficiently large. We note that \eqref{assump1} implies $N\Sinf = O(Ns) = O(N^{1-1/K})$. Therefore, we obtain the simple approximation
\begin{equation}
	\Sinf V_1^\infty\approx \theta\bigl(1-\de_2 s\a\bigr)\,,
\end{equation}
because the term $\frac{1}{\theta Ns\a}$ is swallowed by $O(N^{-K_1/K})$. The asymptotic per-generation response of the mean phenotype is then $\De\bar G^\infty \approx \Th\theta s\a^2\bigl(1-\de_2 s\a\bigr)$; see Corollary 4.8 and Remark 4.9 in GB2024. 

It is of interest to note that the scaling assumption \eqref{assump1} is equivalent to the assumption of `moderately strong selection' used by \citet{Boenkost_etal2021_strong} in their proof of Haldane's approximation for Cannings models, which include the classical Wright-Fisher model. For an illustration of the scaling assumption \eqref{assump1}, we choose $K_1=3/2$ and $K=3$. Then $s=O(N^{-1/3})$, $Ns=O(N^{2/3})$ (which is in contrast to the diffusion approximation), and the error term in \eqref{VGinfty} is $O(N^{-1/2})$.

The above assumptions are needed to derive the error term in \eqref{VGinfty}. 
From \eqref{VGtau} we obtain
\begin{equation}
	V_G^\infty/(\Th\a^2) = \int_0^\infty S^{([t])} \ga([t])\, dt  =  \Sinf \int_0^\infty \ga^\infty(t)\, dt  + D_V \,,
\end{equation}
where $D_V =  \int_0^\infty I_V(t)\, dt$ and $I_V(t) = S^{([t])} \ga(t) - \Sinf \ga^\infty(t)$.
We use the decomposition
\begin{linenomath}\begin{equation}\label{D_V_decomp}
	\int_0^\infty I_V(t)\, dt = \int_0^{T(\ep)} I_V(t)\, dt + \int_{T(\ep)}^\infty I_V(t) \, dt\,,
\end{equation}\end{linenomath}
where $T(\ep)$ is defined in \eqref{Tep_gen_def} and studied in Sect.~\ref{subsect:Tep}. The key points are that (i) if $t\ge T(\ep)$, then $S^{([t])}/\Sinf \le 1+\ep$ and $\int_{T(\ep)}^\infty I_V(t) \, dt$ can be shown to be of order $\ep$ (see eq.~(D.33) in GB2024), and (ii) $T(\ep)$ is a relatively short time, such that the variance contributed by the mutant up to $T(\ep)$ is small, i.e., $\int_0^{T(\ep)} I_V(t)\, dt$ is even smaller (see the derivation of inequality (D.36) in GB2024, which uses the explicit form of the bound \eqref{SnFLgam} to calculate the integral). For the estimates of both integrals in \eqref{D_V_decomp}, the scaling assumption \eqref{assump1} is crucial. Indeed, the proofs show that the choice $\ep = s^{K_1}$, whence $s=O(N^{-K_1/K})$, yields the error term in \eqref{VGinfty}.

The proof of the explicit approximation for the time course of $V_G(\tau)$ in Proposition 4.11 of GB2024 is also based on an analogous time-scale separation and the corresponding inequalities resulting from \eqref{SnFLgam}.

{\bf Acknowledgments.} Thoughtful comments by two anonymous reviewers on a previous version are gratefully acknowledged.

\FloatBarrier

\appendix\appendixpage
\renewcommand{\thesection}{\Alph{section}}

We recommend to consult the supplementary \emph{Mathematica} notebook to check the complicated algebraic computations in the proofs below. For Appendices \ref{sec:Proof_Bin},  \ref{sec:Proof_NB}, and  \ref{sec:Proof_F3}, these are the notebook sections 3, 4, and 5, respectively.

\section{Proof of Theorem \ref{thm:Bin_ge_FL} for the binomial distribution}\label{sec:Proof_Bin}
The proof requires some preparation. Recall that we assume $0\le x \le 1$, $0<\xi<1$, and $\xi^n=\PinfBin$.
We define
\begin{equation}
	\xit = \sum_{k=0}^{n-1} \xi^k = \frac{1-\xi^n}{1-\xi} \; \text{ and }\;  v(x) = \frac{x-\xi^n}{\xi \xit}\,.
\end{equation}
We note that
\begin{equation}\label{inequality_xit}
	\xit > n\xi^{(n-1)/2} 
\end{equation}
(e.g., by the inequality of the arithmetic and geometric means).
Then we obtain
\begin{subequations}\label{v(x)_bounds}
\begin{equation}
	0 < v(x) < \frac{1-\xi^n}{\xi\xit } \; \text{ if }  \xi^n < x < 1
 \end{equation}
and
\begin{equation}\label{v(x)_bounds_b}
	0 > v(x) > v(0) = -\frac{\xi^{n-1}}{\xit}  \ge - \frac{1}{n}\xi^{(n-1)/2} >-\frac{1}{n} \, \text{ if } 0 < x < \xi^n\,.  
\end{equation}
\end{subequations}
Using $v(x)$, we can express the pgfs $\varphiBin$ and $\varphiFL$ as follows:
\begin{equation}\label{varphiBin_redef}
	\varphiBin \Bigl(x;n,\frac{1-\xi}{1-\xi^n}\Bigr) = \xi^n \bigl(1 +v(x)\bigr)^n \,,
\end{equation}
and
\begin{equation}\label{varphiFL_redef}
	\varphiFL(x;\pBin,\rBin) = \xi^n \frac{(1-\xi)(1-\xi^n)+v(x)[n(1-\xi)-\xi(1-\xi^n)]}{(1-\xi)(1-\xi^n)+v(x)[n \xi^n (1-\xi)-\xi(1-\xi^n)]}\,.
\end{equation}

\begin{proof}[Proof of Theorem \ref{thm:Bin_ge_FL}] 
From \eqref{varphiBin_redef} and \eqref{varphiFL_redef} we infer that $\varphiFL(x)\le\varphiBin(x)$ for every $x\in[0,1]$ if and only if
\begin{align}\label{fBinFL}
	\fBinFL(x) &:= \bigl(1 + v(x)\bigr)^n\Bigl((1-\xi)(1-\xi^n)+v(x)[n \xi^n (1-\xi)-\xi(1-\xi^n)]\Bigr) \notag \\
		&\qquad-\Bigl ((1-\xi)(1-\xi^n)+v(x)[n(1-\xi)-\xi(1-\xi^n)]\Bigr) \ge 0 
\end{align}
for every $x\in[0,1]$. In the following, we omit the dependence on $x$ and write $\fBinFL(v)$. By simple rearrangement of \eqref{fBinFL}, we obtain
\begin{align}\label{fBinFL(v)}
	\fBinFL(v) &= \Bigl[\bigl((1 + v)^n - 1\bigr)(1-(1+v)\xi)(1-\xi^n) - n v (1-\xi)\bigl(1-(1 + v)^n\xi^n\bigr)\Bigr] \notag \\
		& = (1-\xi)^2(1-(1+v)\xi) \ffBinFL(v)\,,
\end{align}
where, by using $\frac{1-u^n}{1-u} = \sum_{k=0}^{n-1} u^k$ twice,
\begin{equation}\label{fBinhat}
	\ffBinFL(v) := \bigl((1 + v)^n - 1\bigr)\frac{1}{1-\xi}\sum_{k=0}^{n-1} \xi^k - n v \sum_{k=0}^{n-1}(1+v)^k \frac{1}{1-\xi} \xi^k \,.
\end{equation}
By expansion of $(1-\xi)^{-1}$ and subsequent reordering of the sums we arrive at
\begin{subequations}\label{fBinhatA}
\begin{align}
	\ffBinFL(v) &= v \sum_{j=0}^{n-1} \binom{n}{j+1}v^j \frac{1}{1-\xi}\sum_{k=0}^{n-1} \xi^k  - v \sum_{j=0}^{n-1}  v^j \Bigl(\sum_{k=j}^{n-1} \binom{k}{j} \frac{\xi^k}{1-\xi} \Bigr) \label{fBinhatA_a} \\
	&= v \sum_{j=0}^{n-1} v^j \binom{n}{j+1} \biggl[ \sum_{k=0}^{n-1} (k+1)\xi^k + n\frac{\xi^n}{1-\xi} \biggr] \label{fBinhatA_b}\\
	&\quad - v \sum_{j=0}^{n-1} v^j n \left[ \sum_{k=0}^{n-j-1} \binom{j+k+1}{k} \xi^{j+k} + \binom{n}{j+1}  \frac{\xi^n}{1-\xi} \right] \label{fBinhatA_c}\,.
\end{align}
\end{subequations}
We note that the last terms in each of the brackets  in \eqref{fBinhatA_b} and \eqref{fBinhatA_c} cancel. By setting $k=l-j$ in \eqref{fBinhatA_c}, using $\binom{l+1}{l-j}=\binom{l+1}{j+1}$, and then returning from $l$ to $k$, we obtain
\begin{equation}
	\ffBinFL(v) = v \sum_{j=0}^{n-1} v^j \sum_{k=0}^{n-1} \xi^k \left[ \binom{n}{j+1}(k+1) - n \binom{k+1}{j+1} \right]\,,
\end{equation}
which further simplifies to
\begin{subequations}\label{fBseries}
\begin{equation}\label{fBseries_a}
	\ffBinFL(v) = v^2 \sum_{j=0}^{n-2} \sum_{k=0}^{n-2} v^j \xi^k c_f(n,j,k) \,,
\end{equation}
where
\begin{equation}\label{fBseries_b}
	c_f(n,j,k) =  \binom{n}{j+2}(k+1) - n \binom{k+1}{j+2}  \,.
\end{equation}
\end{subequations}
The coefficients $c_f(n,j,k)$ of $v^j\xi^k$ are always nonnegative. It is sufficient to consider $j+1 \le k \le n-2$. Then
\begin{equation}\label{ratio_product}
	\frac{\binom{n}{j+2}(k+1)}{n \binom{k+1}{j+2}} = \frac{\prod_{i=0}^j (n-1-i)}{\prod_{i=0}^j (k-i)} \ge 1 \,.
\end{equation}
This proves that $\ffBinFL(v)\ge0$ if $v\ge0$, whence $\fBinFL(x)\ge0$ follows if $\xi^n\le x \le 1$.

Now we assume $0<x<\xi^n$, i.e., $v(x)<0$. In the representation \eqref{fBseries} of $\ffBinFL(v)$ we consider two consecutive terms, starting with $j=0$, i.e.,
\begin{equation}
	 v^j \left(c_f(n,j,k) + v c_f(n,j+1,k)\right)\,,
\end{equation}
where $j=0,2,4,\ldots$.
We will prove that this is always positive. Simple calculations show that $\binom{n}{j+3}(k+1) = \binom{n}{j+2}(k+1)\frac{n-j-2}{j+3}$ and $n \binom{k+1}{j+3} = n \binom{k+1}{j+2}\frac{k-j-1}{j+3}$.
Because $v>-\frac{1}{n}$ by \eqref{v(x)_bounds_b}, we obtain
\begin{align}
	&c_f(n,j,k) + v c_f(n,j+1,k) \ge c_f(n,j,k) - \frac{1}{n} c_f(n,j+1,k) \notag \\
		&\qquad= \binom{n}{j+2}(k+1) \left(1- \frac{1}{n}\frac{n-j-2}{j+3} \right) - n \binom{k+1}{j+2} \left(1-\frac{1}{n}\frac{k-j-1}{j+3}  \right) \,.
\end{align}
By \eqref{ratio_product}, the ratio of these two terms simplifies to
\begin{equation}
	\frac{\prod_{i=0}^j (n-1-i)}{\prod_{i=0}^j (k-i)}\cdot \frac{n(j+3)-n-j-2}{n(j+3)-k-j-1}\,,
\end{equation}
which is greater or equal than 
\begin{equation}
	\frac{n-1}{k}\cdot\frac{n(j+3)-n-j-2}{n(j+3)-k-j-1} > 1 \,,
\end{equation}
where the last inequality follows from $(n-1)\bigl(n(j+3)-n-j-2\bigr)-k\bigl(n(j+3)-k-j-1\bigr) = (n-k-1)\bigl((n+1)(j+2)-k\bigr)> 0$ because $k\le n-2$. If the number of such pairs is odd, then $n$ is even. Then the remaining term is positive because $v^n>0$.
This finishes the proof of Theorem \ref{thm:Bin_ge_FL}.
\end{proof}

\begin{remark}
By \eqref{mphi*gaphi<1}, it follows from Theorem \ref{thm:Bin_ge_FL} that $\mBin\gaBin<1$. Here is a simple direct proof. From \eqref{gaBin_xi} we obtain $\mBin\gaBin = \mBin^2 \xi^{n-1}$ and, by \eqref{inequality_xit},
\begin{equation}
	\mBin\gaBin = \mBin^2 \xi^{n-1} = \Bigl(\frac{n}{\xit} \Bigr)^2 \xi^{n-1} < \Bigl(\xi^{-(n-1)/2} \Bigr)^2 \xi^{n-1} =1 \,.
\end{equation}
\end{remark}

\section{Proof of Theorem \ref{thm:NB_ge_FL} for the negative binomial distribution}\label{sec:Proof_NB}
We start by recalling that $\ze^r=\PinfNB \in(0,1)$ and defining
\begin{linenomath}\begin{equation}\label{zet_y_defs}
	\zet = \sum_{k=0}^{r-1} \ze^k \; \text{ and }\; \yx = \frac{\ze^r-x}{\zet} \,. 
\end{equation}\end{linenomath}
Using \eqref{inequality_xit}, we observe that
\begin{equation}\label{y_bounds}
	y(0)=\frac{\ze^r}{\zet} < \frac{1}{r}\ze^{(r+1)/2}<\frac{1}{r} \,, \; y(\ze^r)=0 \,, \; y(1) = \frac{\ze^r-1}{\zet} = -(1-\ze) \,.
\end{equation}
With these abbreviations, we can express $\varphiNB$ in \eqref{phiNB_ze} and $\varphiFL$ in \eqref{frac_lin_gen}, where $\pNB$ and $\rNB$ are given in \eqref{pNB_rNB}, as follows:
\begin{linenomath}\begin{equation}
 \varphiNB\Bigl(x;r,\frac{\ze(1-\ze^r)}{1-\ze^{r+1}}\Bigr) = \ze^r\,\left(\frac{1}{1+\yx} \right)^r
\end{equation}\end{linenomath}
and
\begin{linenomath}\begin{equation}
	\varphiFL(x;\pNB,\rNB) = \ze^r\,\frac{1-x - r\yx}{1-x - r\ze^r \yx}\,.
\end{equation}\end{linenomath}
We note that numerator and denominator are always positive.
These, as all other nontrivial formulas here, are easily verified using the code in Sect.~4 of the supplementary {\it Mathematica} notebook.

\begin{proof}[Proof of Theorem \ref{thm:NB_ge_FL}]
We define
\begin{equation}\label{def_fNB}
	\fNB(x;r,\ze) := \ze^r \left(\varphiFL(x;\pNB,\rNB)^{-1} -  \varphiNB\Bigl(x;r,\frac{\ze(1-\ze^r)}{1-\ze^{r+1}}\Bigr)^{-1} \right)
\end{equation}
and
\begin{equation}\label{def_gNB}
	\gNB(y(x);r,\ze):=\fNB(x;r,\ze)\frac{1-x-r y(x)}{(1-\ze)^2} \,,
\end{equation}
where $1-x-r y(x)>0$. Proving the inequality in \eqref{varphiNB_ge_varphiFL} is equivalent to showing that $\gNB(y(x);r,\ze)>0$ for every $x\in[0,\ze^r)$, $r\ge2$, and $0<\ze<1$. 
In the following we write $y=y(x)$. Using the transformation $x= \ze^r - y\frac{1-\ze^r}{1-\ze}$, we obtain 
\begin{equation}
	\gNB(y;r,\ze) = \fNB\Bigl(\ze^r - y\frac{1-\ze^r}{1-\ze};r,\ze\Bigr) \frac{(1-\ze)(1-\ze^r)+y(1-\ze^r-r(1-\ze))}{(1-\ze)^3}\,,
\end{equation}
which we can rewrite as
\begin{equation} 
	\gNB(y) 
	 = \frac{y\bigl((1+y)^r-1\bigr)\bigl(r(1-\ze)-(1-\ze^r)\bigr) - (1-\ze)(1-\ze^r) \Bigl(\bigl((1+y)^r-1\bigr) - ry \Bigr)}{(1-\ze)^3}
\end{equation}

By the binomial expansion $(1+y)^r-1 = y\sum_{j=0}^{r-1}\binom{r}{j+1}y^j$ and after collection of coefficients of $y^j$, we obtain
\begin{subequations}
\begin{align}
	\gNB(y) &= \frac{y^2}{(1-\ze)^2} \sum_{j=0}^{r-1} y^j\left( \binom{r}{j+1}\Bigl( r-\frac{1-\ze^r}{1-\ze}\Bigr) - \binom{r}{j+2}(1-\ze^r) \right) \\
	&= y^2\sum_{j=0}^{r-1} y^j \binom{r}{j+1}\frac{1}{(1-\ze)^2}\left(\Bigl( r-\frac{1-\ze^r}{1-\ze}\Bigr) - \frac{r-j-1}{j+2}(1-\ze^r) \right) \,.
\end{align}
\end{subequations}
Finally, expansion in terms of $\ze$ yields after appropriate rearrangement
\begin{equation}\label{gNB_final}
	\gNB(y;r,\ze) = y^2 \sum_{j=0}^{r-1} y^j \binom{r}{j+1} c_g(r,j,\ze) \,,
\end{equation}
where
\begin{equation}
	c_g(r,j,\ze) = \frac{1}{2(j+2)}\left( \sum_{k=0}^{r-2} \ze^k(k+1)[2r(1+j)-(2+j)k-2] + \frac{\ze^{r-1}}{1-\ze}r(r+1)j  \right)\,.
\end{equation}
Because $2r(1+j)-(2+j)k-2 \ge 2 + j(2r-k) \ge 2$, we obtain $c_g(r,j,\ze) > 0$ for every $0\le j\le r-1$ and every $\ze\in(0,1)$ if $y>0$. Therefore, $\gNB(y;r,\ze) >0$ if $y>0$ and $\fNB(x;r,\ze)>0$ if $0\le x < \ze^r$.
\end{proof}

We note that the structure of $\gNB$ in \eqref{gNB_final} and the positivity of the coefficients implies that $\gNB$ is strictly convex at $y=0$. Therefore, $\fNB(x)$ is also positive if $x$ is slightly larger than $\ze^r =\PinfNB$.

In the following we show that $\fNB\ge0$ for every $x\in[0,1]$ if $r=2,\ldots,5$. We use the transformation $y=u-(1-\ze)$. Then, by \eqref{y_bounds}, $0\le u \le 1-\ze$ if $\PinfNB = \ze^r \le x \le 1$. 
We use series expansion of $\tilde g_{\rm NB}(u;r,\ze) := \frac{1-\ze}{(u-(1-\ze))^2}\gNB(u-(1-\ze);r,\ze)$. Analytically, this is quite cumbersome due to the complicated structure of $\gNB$. However, {\it Mathematica} performs this task expeditiously (also for $r>5$):
\begin{subequations}
\begin{align}
	\tilde g_{\rm NB}(u;2,\ze) &= u\,, \\
	\tilde g_{\rm NB}(u;3,\ze) &= u(1+4\ze+\ze^2) + u^2(2+\ze)\,, \\
	\tilde g_{\rm NB}(u;4,\ze) &= u(1+4\ze+10\ze^2+4\ze^3+\ze^4) + u^2(2+10\ze+6\ze^2+2\ze^3) \notag \\
		&\quad + u^3(3+2\ze+\ze^2)\,, \\
	\tilde g_{\rm NB}(u;5,\ze) &=  u(1+4\ze+10\ze^2+20\ze^3+10\ze^4+4\ze^5+\ze^6) \notag\\
		&\quad + u^2(2+10\ze+30\ze^2+20\ze^3+10\ze^4+3\ze^5) \notag \\ 
		&\quad + u^3(3+18\ze+13\ze^2+8\ze^3+3\ze^4) + u^4(4+3\ze+2\ze^2+\ze^3)\;.
\end{align}
\end{subequations}
In general, the coefficient of $u$ is $\frac{1}{(1-\ze)^2}(\zet^2-r^2\ze^{r-1})>0$, where the inequality follows from \eqref{inequality_xit}. The coefficient of $u^{r-1}$ is $\frac{1}{1-\ze}(r-\zet)>0$ if $0<\ze<1$.

\begin{remark}
From \eqref{gaNBzet} we obtain $\gaNB=\ze^{r+1}\mNB$, where $\mNB = r/(\ze\zet)$. Therefore, using again \eqref{inequality_xit}, we arrive at
\begin{equation}
	\mNB\gaNB = \mNB^2\ze^{r+1} = \frac{r^2\ze^{r-1}}{\zet^2} < \frac{r^2\ze^{r-1}}{ r^2 \ze^{r-1}} = 1\,.
\end{equation}
\end{remark}

\section{Proofs for distributions with $p_k=0$ for $k\ge4$}\label{sec:Proof_F3}
\begin{proof}[Proof of Lemma \ref{lem:casesF3}]
We recall from Remark \ref{rem:props_p0+_etc}(b) that $\pnr$ is the only (potentially admissible) solution of $p_0=\rFth$. By Remark \ref{rem_props p*}(b), this is the case if and only if \eqref{pnplus_pos} and \eqref{pnr<m} are satisfied. Furthermore, by straightforward algebra,
\begin{linenomath}\begin{equation}\label{p0<rfTh_pnr}
	p_0 < \rFth \;\text{ if and only if }\; p_0 < \pnr \,.
\end{equation}\end{linenomath}

It is not difficult to show that we can write 
\begin{linenomath}\begin{align}\label{gaFth*mFth-1}
	\gaFth \mFth-1 &= \frac{1}{2p_3^2} \Bigl(p_2+3p_3-\sqrt{4p_0p_3+(p_2+p_3)^2}\Bigr)^2 \notag \\
		&\quad\times\biggl(\,\frac{1}{2}-\frac{4p_0p_3+(p_2+p_3)^2+(p_2+3p_3)\sqrt{4p_0p_3+(p_2+p_3)^2}}{4p_3}  \,\biggr)\,.
\end{align}\end{linenomath}
(One way to derive this is to solve $z=\sqrt{4p_0p_3+(p_2+p_3)^2}$ for $p_0$, substitute this expression for $p_0$ in $\gaFth\mFth-1$, and then factorize the resulting quartic polynomial in $z$; see Sect.~5.2 in the {\it Mathematica} notebook). We consider $\gaFth \mFth-1$ as a function of $p_0$. It is well defined if $p_0\ge -(p_2+p_3)^2/(4p_3)$. At this value, it is positive. For sufficiently large $p_0$, it becomes negative and tends to $-\infty$ as $p_0\to\infty$. 
The equation $\gaFth \mFth-1 = 0$ has the solutions $p_0=p_2+2p_3$ (which has multiplicity two) and $p_0=\pnga$. The only potentially admissible solution is $p_0=\pnga$ because $p_0=p_2+2p_3$ yields $\mFth=\gaFth =1$ By Remark \ref{rem_props p*}(c) and (d), $p_0=\pnga$ is admissible if and only if \eqref{pnga>0} and \eqref{pnplus>pnga} hold. By \eqref{p0<rfTh_pnr} and \eqref{pnga<pnr}, this establishes Case (4).

The solution $p_0=p_2+2p_3$ is a critical point of $\gaFth \mFth -1 $, and it is a local maximum if and only if \eqref{pnplus>pnga} is satisfied, i.e., if $\pnga < p_2+2p_3$. In this case, $\gaFth \mFth>1$ if $0<p_0<\pnga$,  and $\gaFth \mFth<1$ if $\pnga< p_0 < p_2+2p_3$. This holds independently of the sign of $\pnga$. Therefore, if $\pnga\le0$, then $\gaFth \mFth<1$ for every admissible $p_0$. If $\pnga > p_2+2p_3$, then $p_0=p_2+2p_3$ is a local minimum of $\gaFth \mFth$, and $\gaFth \mFth > 1$ holds for every admissible $p_0$. We conclude that
\begin{linenomath}\begin{subequations}\label{ga*m><1}
\begin{alignat}{2}
		&\gaFth \mFth > 1 \quad &&\text{if}\quad 0 < p_0 < \min\{\pnga, p_2+2p_3\} \,,\label{ga*m>1} \\
		&\gaFth \mFth < 1 \quad &&\text{if}\quad \max\{\pnga,0\} < p_0 < p_2+2p_3\,. \label{ga*m<1}
\end{alignat}
\end{subequations}\end{linenomath}

Therefore, Case (1) holds by \eqref{p0<rfTh_pnr}, \eqref{pnga<pnr}, and \eqref{ga*m<1}. Case (2) follows from \eqref{pnga<pnr} and \eqref{ga*m<1} because $p_0=\pnr$ is the only (potentially admissible) solution of $p_0=\rFth$. Case (3) follows again from \eqref{p0<rfTh_pnr}, \eqref{pnga<pnr}, and \eqref{ga*m<1}. We deal with the subcases below. Case (4) was already settled above, and Case (5) follows from \eqref{p0<rfTh_pnr}, \eqref{pnga<pnr}, and \eqref{ga*m>1}. Finally, by \eqref{p0<rfTh_pnr} and \eqref{pnga<pnr}, $p_0\ge\rFth$ implies $p_0>\pnr>\pnga$, which is incompatible with $\gaFth \mFth > 1$ by \eqref{ga*m>1}.

Finally, we settle the subcases in Case (3). Indeed, from Remark \ref{rem_props p*}(b), (c), and (d) we already know that $\pnplus<\pnr$ holds if and only if $0< \pnplus$, $\pnga<\pnr$ holds always, and $\pnga<p_2+2p_3$ is equivalent to $\pnga<\pnplus<p_2+2p_3$.
\end{proof}

We start with the proof of Theorem \ref{thm:main_result_F3} after some additional preparation.
We note that the fourth (and all higher) derivatives of $\fFth(x)$ are negative on $[0,1]$. Therefore, $\fFth'''(x)$ is decreasing on $[0,1]$. 
Using the substitution $p_0\to \PinfFth(p_2+p_3+p_3\PinfFth)$ (obtained from \eqref{Pinfty_F3}), we can write
\begin{linenomath}\begin{equation}\label{f(x)_pex3}
	\fFth(x) = \frac{(1-x)(\PinfFth-x)^2\Bigl(-p_3+\bigl(p_2+p_3+2p_3\PinfFth)\bigl(p_2+p_3+p_3\PinfFth+p_3x\bigr)\Bigr)}{1 + (p_2+p_3+2p_3\PinfFth)(\PinfFth-x)} \,.
\end{equation}\end{linenomath}
The denominator is positive on $[0,1]$ because 
\begin{align*}
	&1 + (p_2+p_3+2p_3\PinfFth)(\PinfFth-x) \ge 1 + (p_2+p_3+2p_3\PinfFth)(\PinfFth-1) \\
	&\qquad \ge \PinfFth(p_2+p_3+p_3\PinfFth) + p_2+p_3 + (p_2+p_3+2p_3\PinfFth)(\PinfFth-1) \\
	&\qquad = \PinfFth(2p_2+3p_3\PinfFth) > 0\,,
\end{align*}
where in the second inequality we used $1 \ge p_0+p_2+p_3$.
The numerator is a polynomial of degree four in $x$ with negative leading coefficient.
This informs us that, in addition to the zeroes $\PinfFth$ and 1, which occur by definition, $\fFth(x)$ has at most one additional zero in $[0,1]$. Recalling from \eqref{properties_f_1} that $\PinfFth$ is a critical point of $\fFth$, it has (at least) multiplicity two. Therefore, $\fFth$ can have at most one additional zero in $[0,1]$.

\begin{remark}\label{f''=0}
(a) We recall from \eqref{pnplus_eq} that critical point $\PinfFth$ is a local minimum if and only if $p_0>\pnplus$. 
Therefore, $\PinfFth$ is a local maximum if and only if $p_0 < \pnplus$, which is possible only if \eqref{pnplus_pos} holds.

(b) We recall from Remark \ref{rem:props_p0+_etc}(a) that $\fFth''(\PinfFth) = 0$ if and only $p_0=\pnplus$. Assume $\fFth''(\PinfFth) = 0$. Then $\PinfFth$ and 1 are the only zeroes of $\fFth$ because $\PinfFth$ has multiplicity three. In addition, $\fFth'''(\PinfFth)=3\sqrt{p_3}(p_2-\sqrt{p_3}+3p_3)$, which is positive if $0<\pnplus<p_2+2p_3$ by Remark \ref{rem_props p*}. Hence, $\fFth(x)$ changes sign from negative to positive as $x$ increases from below $\PinfFth$ to above $\PinfFth$. By the considerations above, $\PinfFth(x)$ cannot have additional zeroes, whence $\PinfFth(x)<0$ on $[0,\PinfFth)$ and $\PinfFth(x)>0$ on $(\PinfFth,1)$. 
\end{remark}

\begin{proof}[Proof of Theorem \ref{thm:main_result_F3}]
We distinguish the following cases and recall that $\fFth(x)$ can have three different zeroes only if $\fFth''(\PinfFth)\ne0$.

Case $\fFth'(1)<0$. This is satisfied if and only if $\gaFth\mFth<1$. There are the following subcases: 

(a) $\fFth(x)$ has only the zeroes $\PinfFth$ and 1 in $[0,1]$. Then $\fFth(x)>0$ on $(\PinfFth,1)$, with a local maximum in this interval. Clearly, $\fFth(x)$ cannot have a local maximum at $\PinfFth$. If $\fFth''(\PinfFth)>0$, then $\fFth$ has a local minimum (of 0) at $\PinfFth$. It follows that $\fFth(x)>0$ on $[0,\PinfFth)$ and, in particular, $\fFth(0)>0$, which is equivalent to $p_0> \rFth$. Hence, this occurs precisely in case (1) of Lemma \ref{lem:casesF3} and is displayed in Fig.~\ref{fig_F3}a. 

If $\fFth''(\PinfFth)=0$, i.e., $p_0=\pnplus$, Remark \ref{f''=0} informs us that $\fFth(x)$ changes sign at $\PinfFth$ and is negative below $\PinfFth$, and positive above. By Remarks \ref{rem_props p*}(b) and (d), $p_0=\pnplus$ can occur only if $p_0<\pnr$ and $\pnga < p_0 <p_2+2p_3$. Therefore, $\fFth''(\PinfFth)=0$ can occur only in case (3) of Lemma \ref{lem:casesF3}.

(b) $\fFth(x)$ has a third zero, $x_1\in[0,\PinfFth)$. Then $\fFth(x)>0$ on $(x_1,\PinfFth)$ and on $(\PinfFth,1)$, with local maxima in each of these intervals and a local minimum of 0 at $\PinfFth$, and $\fFth(x)<0$ on $[0,x_1)$. If $x_1=0$, then $\fFth(0)=0$ and $p_0 = \rFth$; see \eqref{pnr_eq}. This is precisely case (2) of Lemma \ref{lem:casesF3} and is displayed in Fig.~\ref{fig_F3}b. 

If $x_1>0$, then $\fFth(0)<0$, i.e., $p_0 < \rFth$ (Fig.~\ref{fig_F3}c). By Lemma \ref{lem:casesF3}, this case applies if \eqref{lem_eq3} holds and, in addition, $\pnga<\pnplus<p_0$. These additional inequalities result from the fact that $\fFth(x)$ has local minimum at $\PinfFth$ (whence $\pnplus<p_0$) and from Remark~\ref{rem_props p*}(d).

(c) $\fFth(x)$ has a third zero, $x_2\in(\PinfFth,1)$. Then $\fFth(x)>0$ on $(x_2,1)$ and $\fFth(x)<0$ on $(\PinfFth,x_2)$. Moreover $\fFth(x)$ must have a local maximum (of 0) at $\PinfFth$, whence $\fFth(x)$ is negative on $[0,\PinfFth)$; see Fig.~\ref{fig_F3}e. In particular, $\fFth(0)<0$, which is equivalent to $p_0 < \rFth$. By Lemma \ref{lem:casesF3}, this case applies if \eqref{lem_eq3} holds and, in addition, $p_0<\pnplus<\pnr$. These additional inequalities result from the fact that $\fFth$ has a maximum at $\PinfFth$ and from Remark~\ref{rem_props p*}(b). 

Case $\fFth'(1)>0$. This is satisfied if and only if $\gaFth\mFth>1$. There are the following subcases:

(a) $\fFth(x)$ has only the zeroes $\PinfFth$ and 1. Because we already know that $\fFth(x)$ cannot change from positive to negative at $\PinfFth$ (when $f''(\PinfFth)=0$), we conclude that $\fFth(x)$ is negative everywhere else and has local maximum of 0 at $\PinfFth$, and a local minimum in $(\PinfFth,1)$; see Fig.~\ref{fig_F3}g. In particular, $\fFth(0)<0$, which is equivalent to $p_0 < \rFth$. This is precisely case (5) of Lemma \ref{lem:casesF3}.

(b) $\fFth(x)$ has a third zero, either in $[0,\PinfFth)$ or in $(\PinfFth,1)$. This is impossible because then  $\fFth(0)\ge0$, i.e., $p_0 \ge \rFth$, which is incompatible with $\gaFth\mFth>1$ by case (6) of Lemma \ref{lem:casesF3}.

Case $\fFth'(1)=0$. This is can be satisfied if only if $\gaFth\mFth=1$, which implies $p_0=\pnga$. Because $p_0<p_2+2p_3$ needs to hold, Remark \ref{rem_props p*} yields $p_0<\pnplus$. Hence $\fFth$ has a local maximum (of 0) at $\PinfFth$. This is precisely case (4) of Lemma \ref{lem:casesF3}. The shape and positivity properties of $\fFth$ are analogous to those in case $\fFth'(1)>0$ (a) (compare Figs.~\ref{fig_F3}f and \ref{fig_F3}g).

The statements about the validity of \eqref{PnFLgam} and \eqref{PnFLgam_reversed} follow because it is sufficient to have $\fFth(x)>0$ ($\fFth(x)<0$) for $0<x < \PinfFth$. The reason is the monotonicity of $\varphi^{(n)}(0)$.
\end{proof}

\singlespacing 
\bibliographystyle{elsarticle-harv}

\end{document}